\documentclass[10pt]{article}
\usepackage{color}
\usepackage{geometry}
\usepackage[algo2e,ruled,lined,titlenumbered]{algorithm2e}

\newcommand{\V}[1]{\underline{#1}}
\newcommand{\M}[1]{\underline{\underline{#1}}}

\newcommand{\lag}{\lambda}

\newcommand{\disp}{\V{u}}

\newcommand{\pos}{\V{x}}
\newcommand{\disptest}{\V{u}^{\star}}

\newcommand{\domega}{\partial \Omega}
\newcommand{\domegau}{\partial \Omega_u}
\newcommand{\domegaf}{\partial \Omega_f}
\newcommand{\ud}{\V{U}_\textrm{D}}
\newcommand{\fd}{\V{f}_\textrm{D}}
\newcommand{\Fd}{\V{F}_\textrm{D}}
\newcommand{\stress}{\M{\sigma}}
\newcommand{\strain}{\M{\epsilon}(\disp)}

\newcommand{\straintest}{\M{\epsilon}(\disptest)}

\newcommand{\dispdiscr}{\V{\mathbf{U}}}
\newcommand{\deltadispdiscr}{\V{\boldsymbol{\delta} \mathbf{U}}}
\newcommand{\Deltadispdiscr}{\V{\boldsymbol{\Delta} \mathbf{U}}}
\newcommand{\resdiscr}{\V{\mathbf{R}}}
\newcommand{\resdiscrred}{\V{\mathbf{R}}_\textrm{R}}

\newcommand{\fint}{\V{\mathbf{F}}_{\textrm{Int}}}
\newcommand{\fext}{\V{\mathbf{F}}_{\textrm{Ext}}}
\newcommand{\zerodiscr}{\V{\mathbf{0}}}

\newcommand{\A}{\M{\mathbf{A}}}
\newcommand{\B}{\V{\mathbf{B}}}
\newcommand{\X}{\V{\mathbf{X}}}
\newcommand{\deltaX}{\V{\boldsymbol{\delta}\mathbf{X}}}
\newcommand{\proj}{\M{\mathbf{P}}}
\newcommand{ \precond}{\M{\widetilde{\mathbf{M}}}}
\newcommand{\optang}{\M{\mathbf{K}}_\textrm{T}}
\newcommand{\optangbar}{\M{\mathbf{\bar{K}}}_\textrm{T}}

\newcommand{\optangred}{\M{\mathbf{K}}_\textrm{T,R}}

\newcommand{\basereduc}{\M{\mathbf{C}}}
\newcommand{\coeffreduc}{\V{\boldsymbol{\alpha}}}
\newcommand{\deltacoeffreduc}{\V{\boldsymbol{\delta \alpha}}}

\newcommand{\schur}{\M{\mathbf{S}}_\mathbf{\textrm{P}}^{(f)}}
\newcommand{\schurforce}{\resdiscr^{(f)}_\textrm{C}}

\usepackage{amssymb} 
\usepackage{amsmath} 
\usepackage{stmaryrd} 
\usepackage{graphicx}
\usepackage{comment}

\usepackage{algorithmic,algorithm}

\title{Local/global model order reduction strategy for the simulation of quasi-brittle fracture
}

\author{P. Kerfriden$^{1}$, J.C. Passieux$^{2}$, S.P.A. Bordas$^{1}$
\\ \\
$^{1}$ Cardiff University, School of Engineering \\
Institute of Mechanics and Advanced Materials \\ 
Queen's Buildings, The Parade, Cardiff CF24 3AA, Wales, United Kingdom \\
$^{2}$ Universit\'e de Toulouse/INSA \\
Institut Cl\'ement Ader (ICA) \\
133 Avenue de Rangueil, F-31077 Toulouse, France
}

\begin{document}
\maketitle

\begin{abstract}
This paper proposes a novel technique to reduce the computational burden associated with the simulation of localised failure. The proposed methodology affords the simulation of damage initiation and propagation whilst concentrating the computational effort where it is most needed, i.e. in the localisation zones. To do so, a local/global technique is devised where the global (slave) problem (far from the zones undergoing severe damage and cracking) is solved for in a reduced space computed by the classical Proper Orthogonal Decomposition, while the local (master) degrees of freedom (associated with the part of the structure where most of the damage is taking place) are fully resolved. Both domains are coupled through a local/global technique. This method circumvents the difficulties associated with model order reduction for the simulation of highly non-linear mechanical failure and offers an alternative or complementary approach to the development of multiscale fracture simulators.

\end{abstract}


\section{Introduction}


Simulating damage initiation and subsequent global structural failure is one of the most active topics in computational mechanics. Several mathematical models and numerical methods have been developed over the years to assess various limit states such as failure due to permanent deformations, cracks or decohesion/delamination, e.g. in composite materials. Yet, these models, be they damage based or relying on discrete cracks are generally computationally  expensive, as they require a fine scale description of the structural and material properties. Therefore, today's engineers are not able to use these state-of-the-art models for routine design. For important recent advances in the treatment of material failure (e.g. discontinuous fracture \cite{moesdolbow1999}, advanced damage models \cite{allixcorigliano1996,jeffersonbennett2009}, damage plasticity models \cite{ladevezelubineau2002,grassljirasek2006} or their combination \cite{rabczukzi2008}, etc.) to become useful in practice, it is thus important to devise techniques which are able to significantly reduce the computational effort required without sacrificing accuracy.

Historically, reducing the computational time associated with solving nonlinear problems in solid mechanics has mainly been addressed by developing homogenisation techniques \cite{fishshek1997,feyelchaboche2000}. In this case, the material properties associated with a material point in a coarse representation of the structure is obtained by averaging of the fine scale material behaviour over a ``representative volume element'' which suppresses the need to resolve the fine scale explicitly over the (possibly large) structure. These methods are widely adopted (see e.g. for a recent description \cite{xiaokarihaloo2009}) and are rather effective in decreasing the computational costs associated with solving fine scale problems. 

However, the application of homogenisation-based methods to failure simulations is not straightforward, as the assumptions required to prove the separation of scales are not satisfied. This may explain why such methods are still under-developed for fracture problems. The ability to coarse-grain fracture information would however lead to a further decrease of computational time by concentrating the computational effort where it is needed most. In the wake of a propagating crack, for example, the full resolution used to resolve fine scale effects close to the crack front is no longer necessary, and a coarse description may be sufficient.  This issue is being intensively investigated by several teams, namely \cite{belytschkoloehnert2007} who developped a coarse graining method called Multiscale Aggregating Discontinuities, and \cite{nguyenlloberas-valls2010,cazescoret2009,verhooselremmers2010} who were able to derive an homogenisation scheme
for lumping the fine scale deffects of Representative Volume Elements (RVEs) into cohesive cracks introduced at the coarse scale.
The solutions proposed so far are however highly problem-dependent and usually do not inherit the solid mathematical foundations from homogenisation theory.

We propose here to follow an alternative route which relies on very recent advances in model order reduction strategies for highly nonlinear structural problems \cite{barraultmaday2004,ryckelynck2005,niroomandialfaro2009,verdonallery2009,kerfridengosselet2010}. We believe that, though rarely studied in the context of failure simulation in nonlinear mechanics, this avenue provides a possibly interesting alternative (or complementary) approach to homogenisation-based solution procedures. This paper is in substance a significant improvement of the work proposed in \cite{kerfridengosselet2010}, dedicated to the fast prediction of the onset of structural failure by adaptive model order reduction. We aim here at preserving the computational efficiency and accuracy obtained in \cite{kerfridengosselet2010} when the size of the damaged zone becomes significant compared to the overall structure. 

The definition of model order reduction can vary depending on authors. We speak here of an approximation of the unknown field in a coarse space spanned by a basis of global vectors (as opposed to locally supported shape functions, in finite elements for example). Model order reduction (MOR) finds its roots in the systems engineering community, where the issue of obtaining a fast and reliable approximation of transfer functions has been extensively addressed over the years. The interested reader can refer to the review proposed in \cite{antoulassorensen2001}. In the context of solid mechanics, most studies have focused on the reduction of dynamics problems. Starting from classical modal basis \cite{rayleighlindsay1877} in a Ritz-Galerkin framework, two main classes of model order reduction methods have been proposed. The first propose a truncation of the modal basis and an approximation of the high-frequency behaviour by static corrections (see for instance \cite{dickensnakagawa1995}). The second use the idea of domain partitioning. Slave degrees of freedom are eliminated by Guyan reduction (or static condensation). The famous Craig-Bampton method \cite{craigbampton1968} is an evolution of this idea, using truncated modal bases to enhance static condensation. These techniques have been extended to the approximation of nonlinear dynamical behaviour (an interesting review is proposed in \cite{kerschenpeeters2009}).

The literature concerning the application of model order reduction to the approximation of the (quasi-)static behaviour of nonlinear structures is a lot sparser. 
Different routes have been investigated to develop such reduced models. The first family of methods consists in taking advantage of the a priori knowledge of a set of representative solutions to fasten the solution process. These representative responses, called snapshots \cite{sirovich1987}, can be for instance a set of solution vectors obtained for given values of structural parameters. The Proper Orthogonal Decomposition (POD) \cite{karhunen1947,loeve1963} can classically be used to extract a relevant global basis from these snapshot vectors, which in turn defines a coarse space where the solution to the structural problem with a different set of parameters can be cheaply searched. 
The second family is not exactly based on reduction by projection. Spectral approaches referred to as Proper Generalized Decomposition (PGD) are a family of solvers capable of building the solution to a multidimensional reference problem as a sum of products of functions of separate variables. The dimensions can be space or time, but also structural parameters. Once this solution is obtained \emph{offline}, the \emph{online} computation simply consists in the evaluation of the separated functions, which are known explicitly. These novel methods are very appealing and were recently applied to linear, stochastic and/or time dependant problems \cite{nouy2007,chinestaammar2010a} (including linear predictions of nonlinear solvers \cite{ladevezepassieux2009}), and to non-linear partial differential equations \cite{nouylemaitre2008,chinestaammar2010a}. Their extension to irreversible constitutive laws such as the ones considered in damage mechanics does not seem to be straightforward, and, as far as we know, has never been adressed before (it is therefore a very interesting possible extension of the present piece of work).

We focus here on the first family of solvers. To the authors' knowledge, the application of model order reduction methods to failure simulations has virtually never been addressed. The main reason for this is that the solution fields in the zones of the structure where damage or failure initiate and propagate are highly parameter-dependent, and cannot be known in advance. Hence, precomputed reduced bases are unable to represent the structural behaviour beyond the onset of strain localisation, damage and fracture. In the particular case of linear elastic fracture mechanics, the authors of \cite{gallandgravouil2010} have proposed a mesh morphing technique that allows to keep an autosimilar description of the singularity during propagation. In this particular case, the use of a snapshot POD is relevant and allow for interesting computational savings. Corrective POD algorithms have also been proposed to alleviate this in the case of plasticity or damage models \cite{ryckelynck2005,kerfridengosselet2010}. Both these techniques are based on global corrections of the initial reduced basis by Krylov algorithms. Though efficient for global nonlinearities, it was shown in \cite{kerfridengosselet2010} that in the case of damage assessment, the number of corrections increases when nonlinearities localise, i.e.: when failure occurs. The purpose of this paper is to circumvent this phenomenon.

The first stage of the proposed method is to only reduce the set of balance equations which exhibits a ``smooth'' nonlinearity. In the context of the simulation of damage, one part of the degrees of freedom, corresponding to highly damaged zones, will be considered as master degrees of freedom, which will be fully solved for, while the remaining will be approximated as a linear combination of reduced basis vectors, which may be precomputed and corrected on the fly. In this respect, the proposed method has strong links with Craig-Bampton methods, or model order reduction by substructure such as the one proposed in \cite{rixen2004,markovicibrahimbegovic2009,ladevezepassieux2009}. A similar approach was also proposed in \cite{ammarpruliere2009}. The second stage of our developments consists in developing an efficient solver for the linearized system of equations obtained by a Newton solution algorithm for the reduced problem. The system is condensed on the master degrees of freedom, and solved by a Krylov algorithm. The global reduced basis is used to efficiently solve the consensed problem: an augmentation technique ensures that the iterations are performed in a space which is orthogonal to the one spanned by the reduced basis. As a result, the augmented Krylov iterations correspond to a two-level solver, the coarse level being the space spanned by the restriction of reduced basis vectors to the fully resolved set of degrees of freedom.

The paper is organized as follows. In the first section, we recall the basics of projection-based model order reduction applied to the solution of nonlinear  quasi-static problems. In section 3, we formally describe the proposed local/global reduction algorithm. This technique is applied to the reduction of damage models in section 4. Finally, we show in the last section that this technique permits to drastically reduce the number of corrections of the reduced basis required to obtain a desired level of accuracy.


\section{Model order reduction by projection: general concepts}


Projection-based model order reduction can be introduced as a means to obtain, at low computational cost, a good approximation to a set of nonlinear equations of the form:
\begin{equation}
\V{\mathbf{F}} \left( \left(\V{\mathbf{X}}_{i} \right)_{i \in \llbracket 0,N \rrbracket} \right) = \zerodiscr
\end{equation}
Where $\V{\mathbf{F}}$ is a vectorial, possibly nonlinear, function of a set of vectorial state variables $\V{\mathbf{X}}$. The expression of this function can take various forms, depending on the application.

\subsection{Problem statement in nonlinear continuum mechanics}

We consider a structure occupying a continuous domain $\Omega$ with boundary $\domega$. It is  subjected to prescribed displacements $\ud$ on its boundary $\domegau$ and to prescribed tractions $\Fd$ on the complementary boundary $\domegaf = \domega \backslash \domegau$, over time interval $[0,T]$. 
Let $\disp$ be the unknown displacement field, which belongs to the space $\mathcal{U}$ of kinematically admissible fields:
\begin{equation}
\label{eq:kinematic}
\mathcal{U}=\left\{\disp\in H^1(\Omega) \ | \  \disp_{| \domegau} = \ud\right\}
\end{equation}
Let $\mathcal{U}^0$ be the associated vector space.
\begin{equation}
\mathcal{U}^0=\left\{\disp\in H^1(\Omega) \ | \  \disp_{| \domegau} = \V{0} \right\}
\end{equation}

Under the assumptions of quasi-static evolution of the structure and small perturbations, the weak form of the equilibrium and constitutive laws read, at any time $t \in [0,T]$:
\begin{equation}
\label{eq:equilibrium}
\begin{array}{l}
\displaystyle \forall \disptest \in \mathcal{U}^0, \ \text{find} \ \disp \in \mathcal{U} \ \text{such that:}
 \\ \displaystyle \qquad \qquad
\int_{\Omega} \stress : \straintest \ d \Omega =  \int_{\Omega} \fd . \disptest \ d \Omega + \int_{\domegaf} \Fd . \disptest \ d \Gamma\\
\displaystyle \qquad \qquad \stress = \stress \left( \left( \strain_{|\tau} \right)_{\tau \leq t} \right)
\end{array}
\end{equation}
where $\stress$ is the Cauchy stress tensor and $\strain$ is the symmetric part of the displacement gradient. 
The constitutive relation between $\stress$ and $\strain$ is nonlinear and described using internal variables (plasticity, damage for instance). It is assumed local and rate-independent.

\subsubsection{Space and time discretization}

Let us perform a standard finite element approximation of the space of unknown displacement field $\mathcal{U}$ (and a similar approximation of the space of test functions $\mathcal{U}^0$):
\begin{equation}
\mathcal{U}^{h}(\Omega) = \left\{ \disp(\V{x}) \ | \ \disp(\V{x}) = \sum_{i=1}^{n_n} N_i(\pos) \, \disp_i  \right\}
\end{equation}
where $n_n$ is the number of nodes, $\pos$ is the position vector, $(N_i)_{i \in \llbracket 1 , n_n \rrbracket}$ are the finite element shape functions associated to each node of the finite element mesh, and $(\disp_i)_{i \in \llbracket 1 , n_n \rrbracket }$ are the nodal values of the displacement field.

The nonlinear solution strategy used to solve the problem over time is a classical time discretization scheme for quasi-static and rate-independent problems.
This procedure consists in finding a set of consecutive solutions at times $(t_n)_{n \in \llbracket 0 , n_t \rrbracket}$.

The introduction of the finite element approximation and time discretization into equation  \eqref{eq:equilibrium} at any time $t_{n+1}$ of the analysis leads to the following nonlinear vectorial equation:
\begin{equation}
\label{eq:discretized_equilibrium}
\fint \left( \Deltadispdiscr , \left(\dispdiscr_{| t_m} \right)_{m \in \llbracket 0 , n \rrbracket} \right) + \fext= \zerodiscr
\end{equation}
where the vector of increment in the nodal unknowns $\Deltadispdiscr \in \mathbb{R}^{n_u}$ ($n_u$ is the number of scalar nodal unknowns) is defined by $\Deltadispdiscr = \dispdiscr_{|t_{n+1}} - \dispdiscr_{|t_{n}}$, $\fint \in \mathbb{R}^{n_u}$ and $\fext \in \mathbb{R}^{n_u}$ are respectively the internal forces resulting from the discretization of the internal virtual work (left-hand side of the first equation of system \eqref{eq:equilibrium}) and the external forces resulting from the discretization of the external virtual work (right-hand side of the first equation of system \eqref{eq:equilibrium}). In the following, the dependency of the internal forces with respect to the history of the unknown fields will not be written explicitly, unless necessary.

The set of successive solution vectors to problem \eqref{eq:discretized_equilibrium} at times $(t_n)_{n \in \llbracket 0 , n_t \rrbracket}$, obtained by a tangent Newton algorithm, will be considered as the reference solution. The algorithms used in the following should provide a solution that is closer to the reference solution, but for a cheaper computational cost. Therefore, the relevance of (i) the initial space and time discretizations and (ii) the chosen Newton solution algorithm will not be discussed in this paper.

\subsection{Reduced solution of the balance equations by projection}

\subsubsection{Principle}

The solution vector is searched for in a space of small dimension (several orders smaller than the number of finite element degrees of freedom). Let us call $\basereduc$ the matrix whose columns form a basis of this space:
\begin{equation}
\basereduc = \left(
\begin{array}{cccc}
\V{\mathbf{C}}^{1} & \V{\mathbf{C}}^{2} & ... & \V{\mathbf{C}}^{n_c}
\end {array}
\right)
\end{equation}
where $n_c$ is the dimension of the reduced space, and $(\V{\mathbf{C}}^{k})_{k \in \llbracket 1 , n_c \rrbracket} \in {(\mathbb{R}^{n_u})}^{n_c}$ are the basis vectors. Applied to the reduction of problem \eqref{eq:discretized_equilibrium}, the increment in the solution field is approximated by:
\begin{equation}
\Deltadispdiscr = \basereduc \, \coeffreduc 
\end{equation}
where we introduced the reduced state variables $\coeffreduc \in \mathbb{R}^{n_c}$. Problem \eqref{eq:discretized_equilibrium} is now overconstrained. One usually tackle this problem by prescribing a Galerkin orthogonality condition: the residual of equation \eqref{eq:discretized_equilibrium} $\resdiscr = \fint ( \dispdiscr_{|t_n} + \basereduc \, \coeffreduc ) + \fext$ is constrained to be orthogonal to any test vector $\deltadispdiscr^{\star} = \basereduc \, \deltacoeffreduc^{\star}$ belonging to the space spanned by the reduced basis vectors $(\V{\mathbf{C}}^{k})_{k \in \llbracket 1 , n_c \rrbracket}$:
\begin{equation}
  \label{eq:ortho_cond}
  \basereduc^T \, \resdiscr = \zerodiscr
\end{equation}
The reduced form of problem \eqref{eq:discretized_equilibrium} is:
\begin{equation}
\label{eq:discr_problem_reduc}
\basereduc^T \left( \fext + \fint( \dispdiscr_{|t_n} + \basereduc \, \coeffreduc ) \right) = \zerodiscr
\end{equation}


\subsubsection{One particular choice of projection space: the proper orthogonal decomposition (POD)}
\label{sec:POD}

Various choices of projection spaces are proposed in the literature. One of the most successful, as far as nonlinear simulations are concerned, is the snapshot-POD \cite{sirovich1987}. This particular technique requires the knowledge of a representative family of solutions to the global problem. 
This set of vectors $(\V{\mathbf{S}}^{k})_{k \in \llbracket 1 , n_s \rrbracket}$ is called \emph{snapshot}. 
The aim is to find an orthonormal basis $(\V{\mathbf{C}}^{k})_{k \in \llbracket 1 , n_c \rrbracket}$, of dimension $n_c$ smaller than $n_s$ such that the distance between spaces $\textrm{Im}(\mathbf{\M{C}})$ and $\textrm{Im}(\M{\mathbf{S}})$ is minimum. This distance can be quantified by the following metric:
\begin{equation}
  \label{eq:fonctionnelle}
 J(\V{\mathbf{C}}^{1},\dots,\V{\mathbf{C}}^{n_c}) = \sum_{j=1}^{n_s} \left\Vert  \V{\mathbf{S}}^{j} - \sum_{i=1}^{n_c} \left( {\V{\mathbf{C}}^{i}}^T \V{\mathbf{S}}^{j} \right) \V{\mathbf{C}}^i \right\Vert^2
\end{equation}
which must be minimized under the constraint of orthonormality of $(\V{\mathbf{C}}^{k})_{k \in \llbracket 1 , n_c \rrbracket}$. Problem \eqref{eq:fonctionnelle} can be reformulated as an unconstrained formulation:
\begin{equation}
  \label{eq:lagfonctionnelle}
 L(\V{\mathbf{C}}^{1},\dots,\V{\mathbf{C}}^{n_c},\lag_{11},\dots,\lag_{n_c n_c}) = J(\V{\mathbf{C}}^{1},\dots,\V{\mathbf{C}}^{n_c}) + \sum_{i=1}^{n_c} \sum_{j=i}^{n_c} \lag_{ij} \left( {\V{\mathbf{C}}^{i}}^T \V{\mathbf{C}}^{j} - \delta_{ij}\right)
\end{equation}
where $(\lag_{ij})_{i,j \in \llbracket 1 , n_c \rrbracket}$ are Lagrange multipliers and $\delta_{ij}$ the Kronecker symbol. 
The $n_c(n_c+1)/2$ optimality conditions for maximizing the Lagrangian with respect to the lagrange multiplier naturally ensure the orthogonality of the basis:
\begin{equation}
  \label{eq:minifonclag}
\frac{\partial L}{\partial \lag_{ij}}=0 \quad \rightarrow  \quad {\V{\mathbf{C}}^{i}}^T \V{\mathbf{C}}^{j} = \delta_{ij}
\end{equation}
Whereas the $n_c$ optimality conditions for minimizing the Lagrangian with respect to the basis vectors reads:
\begin{equation}
  \label{eq:minifonc2}
\frac{\partial L}{\partial \V{\mathbf{C}}^{i}}=\zerodiscr \quad \rightarrow  \quad \sum_{j=1}^{n_s} \V{\mathbf{S}}^{j} \left( {\V{\mathbf{C}}^{i}}^T \V{\mathbf{S}}^{j} \right) = \lag_{ii} \; \V{\mathbf{C}}^i \qquad \text{and} \quad \lag_{ij}=0 \quad \text{for} \quad i\neq j
\end{equation}
Defining $\lag_i=\lag_{ii}$ and $\M{\mathbf{S}} = \left( \begin{array}{cccc} \V{\mathbf{S}}^1 & \V{\mathbf{S}}^2 & ... &  \V{\mathbf{S}}^{n_s} \end{array} \right) \in \mathbb{R}^{n_u} \times \mathbb{R}^{n_s}$, equations \eqref{eq:minifonc2} can be written as the following $n_u\times n_u$ eigenvalue problem:
 \begin{equation}
  \label{eq:minifonc}
\M{\mathbf{S}} \, \M{\mathbf{S}}^T \V{\mathbf{C}}^{i} = \lag_{i} \; \V{\mathbf{C}}^i 
\end{equation}
This problem is equivalent to computing the Singular Value Decomposition (SVD) 
of $\M{\mathbf{S}}$, where the singular values $s_i$ are such that $s_i^2=\lag_i \geq 0$. If $n_s \ll n_u$ (which is the case in general when applying the snapshot-POD),  replacing problem \eqref{eq:minifonc} by the eigenvalue problem on the $n_s\times n_s$ operator $\M{\mathbf{S}}^T \M{\mathbf{S}}$ is computationally cheaper, and reads:
 \begin{equation}
  \label{eq:pbvpplusfacile}
 \V{\mathbf{C}}^{i} = \lag_{i}^{-\frac{1}{2}} \, \M{\mathbf{S}} \, \V{\mathbf{V}}^i \qquad \text{with} \qquad \M{\mathbf{S}}^T \M{\mathbf{S}} \, \V{\mathbf{V}}^{i} = \lag_{i} \; \V{\mathbf{V}}^i
\end{equation}
The error associated with the truncation of the SVD at rank $n_c$ is quantifyed by the distance $J(\V{\mathbf C}^1,...,\V{\mathbf C}^{n_c})$ which can be shown to be equal to the sum of the truncated eigenvalues:
\begin{equation}
J(\V{\mathbf C}^1,...,\V{\mathbf C}^{n_c})=\sum_{i=n_c+1}^{n_s} \lag_i
\end{equation}
This result provides a precious tool for choosing the most relevant basis. Indeed, only the basis vectors $\V{\mathbf C}_i$ associated with the largest $\lag_i$ (in practice, those such that $\lag_i/\lag_{max}>\varepsilon$, $\varepsilon$ being a small parameter depending on the required accuracy) are selected as reduced basis vectors. This ensures that for an given number $n_c$ of basis vectors, the approximation of the snapshot space is maximised in the sense of \eqref{eq:fonctionnelle}.

\subsection{Nonlinear solution algorithm}

At any time step, problem \eqref{eq:discr_problem_reduc} can be solved by classical Newton algorithms (successive linearizations). Iteration $i+1$ of the tangent Newton method consists in solving the following linear predictor:
\begin{equation}
\label{eq:prediction_NR_red}
\optangred^i \, \deltacoeffreduc^{i+1} = -\resdiscrred^{i}
\end{equation}
where the increment in the reduced state variables is $\deltacoeffreduc^{i+1} =\coeffreduc ^{i+1} - \coeffreduc^{i}$, and the residual and tangent operator obtained from the knowledge of the solution at iteration $i$ read (see \cite{kerfridengosselet2010} for more details):
\begin{equation}
\left\{ \begin{array}{l}
\displaystyle \resdiscrred^{i} = \basereduc^T \left( \fint(\dispdiscr_{| t_n}+\basereduc \, \coeffreduc^{i}) +\fext \right) \\
\displaystyle {\optangred^{i}} = \basereduc^T \, \left. \frac{\partial \fint(\dispdiscr_{| t_n}+\basereduc \, \coeffreduc)}{\partial \coeffreduc} \right|_{ \coeffreduc = \coeffreduc^{i}}
\end{array} \right.
\end{equation}


\section{Local/global model order reduction strategy}
\label{sec:LocGlob}


Projection-based model order reduction techniques introduce a global approximation of the displacement. As such, even adaptive versions \cite{ryckelynckchinesta2006,kerfridengosselet2010} of these methods are not well suited to the analysis of localised nonlinearities. We propose in the following to use reduction techniques by projection for the ``weakly nonlinear" equations of reference problem \eqref{eq:discretized_equilibrium}, namely equations for which reduced order modelling (ROM) is relevant, while the remaining equations will be solved directly.

\subsection{Displacement approximation}

The principle of the proposed local/global strategy is to split the unknown solution vector into two parts, only one of them being approximated as a linear combination of reduced basis vectors. We shall use superscript $^{(r)}$ for the approximated ``slave'' part of the solution vector, while superscript $^{(f)}$ (which stands for ``fully resolved") will be used to denote its complementary ``master'' part.

More precisely, let $\widetilde{\Deltadispdiscr}$ be the unknown increment vector at time $t_{n+1}$, the numbering being reorganised as follows:
\begin{equation}
\widetilde{\Deltadispdiscr} = 
\left( \begin{array}{c}
\widetilde{\Deltadispdiscr}^{(r)}  \\
\widetilde{\Deltadispdiscr}^{(f)} 
\end{array} \right)
=
\left( \begin{array}{c}
\M{\mathbf{E}}^{(r)} \\
\M{\mathbf{E}}^{(f)}
\end{array} \right)  \Deltadispdiscr 
\end{equation}
In the above notations, $\widetilde{\Deltadispdiscr}^{(r)} \in \mathbb{R}^{n_r}$ and $\widetilde{\Deltadispdiscr}^{(f)} \in \mathbb{R}^{n_f}$, with $n_r + n_f = n_u$ and $n_f$ much smaller than $n_r$. $\M{\mathbf{E}}^{(r)} \in  \{ 0 , 1 \}^{n_r} \times \{ 0 , 1 \}^{n_u}$ and $\M{\mathbf{E}}^{(f)} \in \{ 0 , 1 \}^{n_f} \times \{ 0 , 1 \}^{n_u}$ are two boolean extractors (rectangular matrices with one 1 per line, the other coefficients being null).

$\widetilde{\Deltadispdiscr}^{(r)}$ is approximated as a linear combination of global vectors $(\V{\mathbf{C}}^{k})_{k \in \llbracket 1 , n_c \rrbracket} \in {(\mathbb{R}^{n_u})}^{n_c}$:
\begin{equation}
\widetilde{\Deltadispdiscr}^{(r)} = \left( {\M{\mathbf{E}}^{(r)}} \basereduc \right) \coeffreduc
\end{equation}

where $\coeffreduc \in \mathbb{R}^{n_c}$ is the vector of reduced degrees of freedom. In the initial numbering, the unknown vector can be expressed in the form:
\begin{equation}
\begin{array}{l}
\displaystyle \Deltadispdiscr = \Deltadispdiscr^{(r)} + \Deltadispdiscr^{(f)}  \\
\displaystyle \textrm{where} 
\left\{ \begin{array}{l}
\displaystyle \Deltadispdiscr^{(r)} = {\M{\mathbf{E}}^{(r)}}^{T} \widetilde{\Deltadispdiscr}^{(r)} \\
\displaystyle \Deltadispdiscr^{(f)} = {\M{\mathbf{E}}^{(f)}}^{T} \widetilde{\Deltadispdiscr}^{(f)}
\end{array} \right. \end{array}
\end{equation}
Introducing projector $\M{\mathbf{P}}^{(r)} ={\M{\mathbf{E}}^{(r)}}^T \M{\mathbf{E}}^{(r)}$, the approximation of the displacement increment finally reads:
\begin{equation}
\Deltadispdiscr = \M{\mathbf{P}}^{(r)} \basereduc \, \coeffreduc + {\M{\mathbf{E}}^{(f)}}^{T} \widetilde{\Deltadispdiscr}^{(f)}
\end{equation}

For the sake of legibility, we define the new state vector $\X$ at time $t_{n+1}$:
\begin{equation}
\X = 
\left( \begin{array}{c}
\coeffreduc \\
\widetilde{\Deltadispdiscr}^{(f)}
\end{array} \right) 
\end{equation}
Using the above definition, the approximation of the displacement increment is now:
\begin{equation}
\Deltadispdiscr = \M{\mathbf{A}}\; \X \qquad \text{where} \qquad \M{\mathbf{A}} = \left( \M{\mathbf{P}}^{(r)}  \basereduc^\textrm \qquad {\M{\mathbf{E}}^{(f)}}^T \right)
\end{equation}

\subsection{Locally reduced set of balance equations}


At that point, we need to define how the overconstrained balance equations \eqref{eq:discretized_equilibrium} will be solved for. Let us restrict our numerical investigations to the Galerkin procedure. The balance equation are therefore required to be orthogonal to any test vector:
\begin{equation}
\deltadispdiscr^\star = \M{\mathbf{A}} \; 
\left( \begin{array}{c}
\coeffreduc^\star \\
\widetilde{\deltadispdiscr}^{(f) \star}
\end{array} \right)
\end{equation}
where ${\deltacoeffreduc}^{\star} \in \mathbb{R}^{n_c}$ and $\widetilde{\deltadispdiscr}^{(f) \star} \in \mathbb{R}^{n_f}$ are arbitrary test state variables.
The orthogonality requirement applied to the initial nonlinear set of equations \eqref{eq:discretized_equilibrium} yields the following nonlinear system:
\begin{equation}
\label{eq:discr_problem_reduc_loc_glo}
\resdiscrred(\X) = 
\M{\mathbf{A}}^T
\left( \fint \left( \Deltadispdiscr (\X) + \dispdiscr_{|t_n} \right) + \fext \right) = \zerodiscr
\end{equation}

\subsection{Newton solution scheme}

At each time step of the time discretization scheme, the reduced problem \eqref{eq:discr_problem_reduc_loc_glo} is solved by a Newton algorithm. The $({i+1})^{th}$ Newton iteration consists in finding $\deltaX^{i+1}$ satisfying the following set of $n_c + n_\textrm{f}$ linearized equation:
\begin{equation}
\label{eq:linearized_reduced}
\displaystyle \left. \frac{\partial \resdiscrred(\X) }{\partial \X} \right|_{ \X = \X^{i} }  \deltaX^{i+1}  = - \resdiscrred^{i}
\end{equation}
where $\deltaX^{i+1}=\X^{i+1}-\X^{i}$ and $\resdiscrred^{i}$ is the residual of problem \eqref{eq:discr_problem_reduc_loc_glo} computed at iteration $i$:
\begin{equation}
\resdiscrred^{i} = \resdiscrred(\X^{i})
\end{equation}

The expression of the tangent operator $\optangred^i = \displaystyle \left.  \frac{ \partial \resdiscrred(\X) }{\partial \X} \right|_{ \X^{i} }$ is obtained by differentiation of  equation \eqref{eq:discr_problem_reduc_loc_glo} with respect to the set of reduced state variables, in any arbitrary direction $\deltaX$ such that:
\begin{equation}
\deltaX =\left( \begin{array}{c} \deltacoeffreduc \\ \widetilde{\deltadispdiscr}^{(f)} \end{array} \right)
\end{equation}
This differentiation reads:
\begin{equation}
\optangred^i \, \deltaX =  \left. \frac{\partial\resdiscrred(\X)}{\partial \X} \right|_{ \X^{i} } \, \deltaX
\end{equation}
Using the tangent operator of the full system of equations \eqref{eq:discretized_equilibrium}, and assuming that $\displaystyle \frac{\partial \fext}{ \partial \Deltadispdiscr}=\zerodiscr$, on can write:
\begin{equation}
\optangred^i \, \deltaX = \M{\mathbf{A}}^T \, \left. \frac{\partial \fint (\Deltadispdiscr) }{\partial \Deltadispdiscr} \right|_{ \Deltadispdiscr = \M{\mathbf{A}} \, \X^{i} } \left.  \frac{\partial \Deltadispdiscr}{\partial \X} \right|_{ \Deltadispdiscr = \M{\mathbf{A}} \, \X^{i} }  \, \deltaX
\end{equation}
Using the notation $\displaystyle \optang^i = \left. \frac{\partial \fint}{\partial \Deltadispdiscr} \right|_{ \Deltadispdiscr (\X^{i}) }$, the previous equation yields the expression of the reduced tangent operator as a function of the tangent of the full system of equations:
\begin{equation}
\optangred^i \, \deltaX = \M{\mathbf{A}}^T \optang^i \M{\mathbf{A}} \, \deltaX
\end{equation}
Introducing the above expression of the tangent into equation \eqref{eq:linearized_reduced}, and expanding the expression of operator $\M{\mathbf{A}}$, one finally obtains the following linearized system:
\begin{equation}
\label{eq:linearized_reduced_matrix}
\begin{array}{l}
\left( \begin{array}{cc}
\displaystyle \basereduc^{T} \M{\mathbf{P}}^{(r)} \optang^i \M{\mathbf{P}}^{(r)} \basereduc & 
\displaystyle \basereduc^{T} \M{\mathbf{P}}^{(r)} \optang^i {\M{\mathbf{E}}^{(f)}}^T
\\
\displaystyle  {\M{\mathbf{E}}^{(f)}} \optang^i \M{\mathbf{P}}^{(r)} \basereduc  &
\displaystyle {\M{\mathbf{E}}^{(f)}} \optang^i  {\M{\mathbf{E}}^{(f)}}^T
\end{array} \right)
\deltaX^{i+1}
 \\ 
\qquad \qquad =
- \left( \begin{array}{c}
\displaystyle  \basereduc^T  \M{\mathbf{P}}^{(r)}  \\
\displaystyle  \M{\mathbf{E}}^{(f)} 
\end{array} \right) \left( \fint \left( \Deltadispdiscr (\X^{i}) + \dispdiscr_{|t_n}  \right) + \fext \right)
\end{array}
\end{equation}

Notice that problem \eqref{eq:linearized_reduced_matrix} is considerably reduced compared to a direct linearisation of \eqref{eq:discretized_equilibrium}. Indeed, $\optangred^i \in \mathbb{R}^{n_c} \times \mathbb{R}^{n_f}$ where $n_c \ll n_u$ and, applied to the analysis of localised nonlinearities, we can reasonably expect that $n_f$ is at least one order of magnitude smaller than $n_u$, which means that only a few of the balance equations are strongly nonlinear. This assumption will be validated in section \ref{sec:application_damage}.

This framework is very general. It does not rely on domain decomposition but on a splitting of the balance equations. Yet, this technique will be particularized in section \ref{sec:application_damage}, and we will show that using the ideas inspired from domain decomposition methods \cite{farhatroux1994,kerfridenallix2009} and local/global methods \cite{gendreallix2009,allixkerfriden2010b} provides efficient splitting in the case of damage mechanics.

\subsection{Two-level solution of the linearized systems}

The successive solution to \eqref{eq:linearized_reduced} could be obtained by a direct solver. This would not affect the results given in the following. However, we believe that, for this particular application, an iterative solver has several advantages:
\begin{itemize}
\item Iterative solvers, such as a conjugate gradient (see \cite{saad2003} for instance), do not require any factorisation. Factorisations  usually deteriorate the sparsity of algebraic systems, which results in an increase in the memory requirements. Newton-Krylov methods do not, in theory \cite{knollkeyes2004}, even require the assembly of the stiffness operator.
\item Direct solvers are less versatile than iterative ones, and are usually a limiting point to extend numerical solution schemes to parallel computing.
\item Solving exactly the successive linearized systems \eqref{eq:linearized_reduced} is not required. This idea has led to intensive work on Inexact Newton Methods \cite{demboeisenstat1982}, but this issue will not be addressed in this paper. The interested reader can refer to \cite{kerfridenallix2009,kerfridengosselet2010} for studies on similar topics in the scope of domain decomposition and model order reduction.
\item More importantly, a basic direct solver such as Cholesky factorizations will not make use of any mechanically relevant information obtained during the previous iterations of the nonlinear solver. In the chosen iterative algorithm, the part of the displacement field that is ``fully resolved'' will be solved for efficiently  using the reduced basis as a preconditioner, following the ideas proposed for instance in \cite{rislerrey2000,knollkeyes2004,kerfridenallix2009}.
\end{itemize}

We use here a Krylov algorithm to efficiently find a solution to linearized system of equation \eqref{eq:linearized_reduced}. Three stages of preconditioning are performed: a condensation, a projection and a crude preconditioning of the resulting system.

\subsubsection{Condensation}

Let us condense the linearized problem \eqref{eq:linearized_reduced_matrix} on the master (i.e.: ``fully resolved'') degrees of fredom:
\begin{equation}
\label{eq:condensed}
\schur \, \widetilde{\deltadispdiscr}^{(f)} = \schurforce
\end{equation}
where the superscript corresponding to the current and previous Newton iterations have been dropped for the sake of legibility. $\schur$ is the primal Schur complement and $\schurforce$ the condensed residual. They are defined by:
\begin{equation}
\left\{ \begin{array}{l}
\displaystyle \schur = \optangred^{(ff)} - \optangred^{(fr)} \left( \optangred^{(rr)} \right)^{-1} \optangred^{(rf)}  \\
\displaystyle \schurforce =  \optangred^{(fr)} \left( \optangred^{(rr)} \right)^{-1} \resdiscrred^{(r)} - \resdiscrred^{(f)}
\end{array} \right.
\end{equation}

Where we have used the usual block notations:
\begin{equation}
\label{eq:linearized_reduced_matrix2}
\optangred = 
\left( \begin{array}{cc}
\displaystyle \optangred^{(rr)} &  \displaystyle \optangred^{(rf)}
\\
\displaystyle \optangred^{(fr)} & \displaystyle \optangred^{(ff)}
\end{array} \right)
\end{equation}

\begin{equation}
\resdiscrred = 
\left( \begin{array}{c}
\resdiscrred^{(r)}
\\
\resdiscrred^{(f)}
\end{array} \right)
\end{equation}

The state variables corresponding to the reduced part of the displacement  vector increment have been eliminated from the system of equations and can be retrieved as follows:
\begin{equation}
\deltacoeffreduc = - \left( \optangred^{(rr)} \right)^{-1} \left( \resdiscrred^{(r)} + \optangred^{(rf)} \, \widetilde{\deltadispdiscr}^{(f)} \right)
\end{equation}

$\left( \optangred^{(rr)} \right)^{-1}$ is an operator of very small size, which can be explicitly computed prior to applying the iterative Krylov solver. However, in the case where a very large number of reduced basis vectors is used to approximate the unknown field in the smooth region, computing the Schur-complement/vector products corresponding to the construction of the successive Krylov basis vectors will become computationally more efficient.

\subsubsection{Iterative solution to the condensed problem}
\label{sec:iterativeCG}

The classical projected (or augmented) conjugate gradient \cite{dostal1988,farhatroux1994} is applied to the approximate solution of \eqref{eq:condensed} (the linearized operator is assumed symmetric, positive and definite). The chosen augmentation space is the space spanned by the reduced basis vectors $\textrm{Im} (\basereduc^{(f)} ) $ with the restriction $ \basereduc^{(f)} =  {\M{\mathbf{E}}^{(f)}}  \basereduc$.

The starting point of augmented Krylov solvers is to decompose the unknown solution increment into two supplementary spaces:
\begin{equation}
\label{eq:uncoupl_CG}
\begin{array}{l}
\displaystyle \widetilde{\deltadispdiscr}^{(f)} = \deltadispdiscr^{(f)}_\textrm{C} + \deltadispdiscr^{(f)}_\textrm{K}
\\
\displaystyle \textrm{where} \quad
\left\{ \begin{array}{l}
\displaystyle \deltadispdiscr^{(f)}_\textrm{C} \in \textrm{Im} \left( \basereduc^{(f)} \right) \\ 
\displaystyle \deltadispdiscr^{(f)}_\textrm{K} \in \textrm{Im} (\basereduc^{(f)})^{\bot_{\schur}} = \textrm{Ker} \left( {\basereduc^{(f)}}^T \schur \right)
\end{array} \right.
\end{array}
\end{equation}
$\bot_{\schur}$ designing the $\schur-$orthogonality, which is ensured by introducing the classical projector:
\begin{equation}
\proj = \M{\mathbf{I}}_\textrm{d} - \basereduc^{(f)} ({\basereduc^{(f)}}^T \schur \basereduc^{(f)}  )^{-1} {\basereduc^{(f)}}^T \schur
\end{equation}
and writing that $\deltadispdiscr^{(f)}_\textrm{K} = \proj \, \deltadispdiscr^{(f)}_\textrm{K}$.

This separation of the search space into two subspaces $\textrm{Im} \left( \basereduc^{(f)} \right)$ and $\textrm{Im}(\proj)$ in direct sum leads to the following uncoupled equations:
\begin{equation}
\label{eq:line_ini}
\begin{array}{l}
\displaystyle 
\deltadispdiscr^{(f)}_\textrm{C}= \basereduc^{(f)} ({\basereduc^{(f)}}^T  \schur \, \basereduc^{(f)})^{-1} {\basereduc^{(f)}}^T \schurforce
\\
\displaystyle 
\left( \proj^T \schur \, \proj \right)  \mathbf{\deltadispdiscr}^{(f)}_\textrm{K} = \proj^T \schurforce
\end{array} \end{equation}
The first line of equation \eqref{eq:line_ini} is a coarse initialization of the projected conjugate gradient. The second line is the linear prediction of the full problem projected on $\textrm{Im} (\basereduc^{(f)})^{\bot_{\schur}}$. This system is symmetric and can be solved by a preconditioned conjugate gradient. Hence one solves iteratively:
\begin{equation}
\precond^{-1} \left( \proj^T \schur  \, \proj \right) \deltadispdiscr^{(f)}_K = \precond^{-1} \, \proj^T  \schurforce
\end{equation}
where  $ \precond^{-1} $ is the left preconditioner (symmetric, definite and positive). In our test cases, $\precond$ is a diagonal matrix whose entries are the elements of the diagonal of $\schur$. Algorithm~\ref{APCG} gives the classical implementation of an augmented preconditioned conjugate gradient.

\begin{algorithm2e}[ht]\caption{Augmented preconditioned conjugate gradient for the solution to the condensed linearized reduced problem}\label{APCG}
Compute 
$\left\{ \begin{array}{l}
\displaystyle \schur = \optangred^{(ff)} - \optangred^{(fr)} \left( \optangred^{(rr)} \right)^{-1} \optangred^{(rf)}  \\
\displaystyle \schurforce = \resdiscrred^{(f)} - \optangred^{(fr)} \left( \optangred^{(rr)} \right)^{-1} \resdiscrred^{(r)} \\
\basereduc^{(f)} =  {\M{\mathbf{E}}^{(f)}}  \basereduc \\
({\basereduc^{(f)}}^T \schur  \basereduc^{(f)})^{-1}
\end{array} \right.$
\\
$\deltadispdiscr^{(f)}_\textrm{C} = \basereduc^{(f)}({\basereduc^{(f)}}^T \schur \basereduc^{(f)})^{-1} {\basereduc^{(f)}}^T \schurforce$ \;
$\V{\mathbf{R}}_{\textrm{CG},0}=\schurforce-\schur \deltadispdiscr^{(f)}_\textrm{C}=\proj^T \schurforce$\;%
$\V{\mathbf{Z}}_0 = \proj \, \precond^{-1} \V{\mathbf{R}}_{\textrm{CG},0} \quad $, $\V{\mathbf{W}}_0=\V{\mathbf{Z}}_0 \quad$ and $\quad  {\mathbf{\deltadispdiscr}^{(f)}_\textrm{K,0}}=0$\;%
\For{$j=1,\ldots,n$}{%
  $\alpha_{j-1}=(\V{\mathbf{R}}_{\textrm{CG},j-1},\V{\mathbf{W}}_{j-1})/(\schur  \V{\mathbf{W}}_{j-1},\V{\mathbf{W}}_{j-1})$  \\
  $ {\mathbf{\deltadispdiscr}^{(f)}_\textrm{K,j}}= {\mathbf{\deltadispdiscr}^{(f)}_\textrm{K,j-1}} +\alpha_{j-1} \V{\mathbf{W}}_{j-1}$\\
  $\V{\mathbf{R}}_{\textrm{CG},j} = \V{\mathbf{R}}_{\textrm{CG},j-1} -\alpha_{j-1} \schur  \V{\mathbf{W}}_{j-1}$\\
  $\V{\mathbf{Z}}_{j} = \proj \, \precond^{-1} \V{\mathbf{R}}_{\textrm{CG},j}$\\
  $\beta_j=(\schur \V{\mathbf{W}}_{j-1},\V{\mathbf{Z}}_{j})/(\V{\mathbf{W}}_{j-1},\schur \V{\mathbf{W}}_{j-1})$ \\
  $\V{\mathbf{W}}_{j}=\V{\mathbf{Z}}_{j} - \beta_j \V{\mathbf{W}}_{j-1}$

}%
\end{algorithm2e}

\subsubsection{Interpretation and comments on the expected additional costs}

\begin{figure}[p]
       \centering
       \includegraphics[width=1. \linewidth]{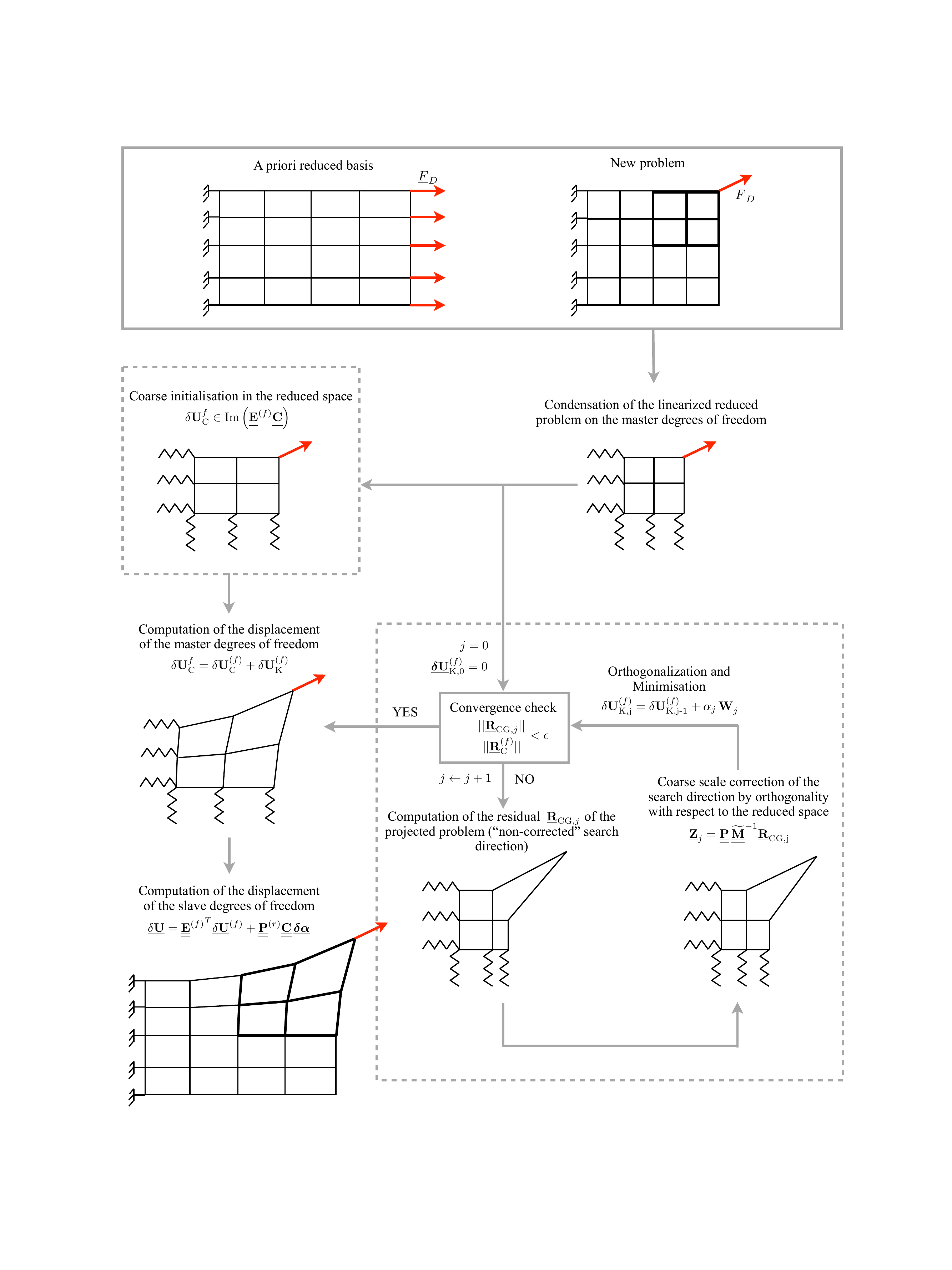}
       \caption{Graphical representation of the two-level solution algorithm.}
       \label{fig:schematic}
\end{figure}

The three preconditioning techniques that we have detailed in the previous sections define a cheap structural preconditioner. A schematic representation of this algorithm is given in figure (\ref{fig:schematic}). \\

From an algebraic point of view, it can easily be proved that the condensation has the same effect on the conjugate gradient iterations as a left block preconditioner for system \eqref{eq:linearized_reduced_matrix}:
\begin{equation}
{\precond_{\textrm{Condens}}}^{-1} = 
\left( \begin{array}{cc}
\displaystyle \left( \optangred^{(rr)} \right)^{-1} & \displaystyle \M{\mathbf{0}} \\
\displaystyle \M{\mathbf{0}} & \displaystyle \M{\mathbf{0}} 
\end{array} \right)
\end{equation}
From a mechanical point of view, condensing the problem on the master degrees of freedom ensures that any solution provided by the iterative solver exactly solves the reduced part of the balance equations.

The augmentation of the conjugate gradient can also be seen as a left preconditioner for the condensed system of equations (second line of system \eqref{eq:line_ini}), which is associated with an initialization of the iterative alogrithm. The initialization $\deltadispdiscr^{(f)}_\textrm{C}$ (first line of system \eqref{eq:line_ini}) is the best solution of the condensed problem in the restriction of the pre-computed reduced basis to the master degrees of freedom (in the sense of a $(\schur)^{-1}$-norm). Taking into account the previous comments on the interpretation of the condensation step, this initialisation is exactly the solution to the linearized system of equation \eqref{eq:prediction_NR_red} on the master degrees of freedom. In other words, the initialization of the iterative algorithm provides the solution of the reduced linearized problem which would have been obtained without local/global enhancement. No significant additional costs compared to the basic projection-based reduced order modelling are involved so far. In the case where the solution which is looked for is correctly approximated in the reduced space (e.g.: ``smooth'' nonlinearity, or unloading of the structure), no conjugate gradient iteration is required.

If conjugate gradient iterations are required, the augmentation acts as a global correction of the successive search directions (see figure (\ref{fig:schematic})). It ensures that any solution of the form \eqref{eq:uncoupl_CG} still satisfies the reduced balance equations without local/global enhancement. Therefore, the complementary solution $\deltadispdiscr^{(f)}_\textrm{K}$ can be interpreted as a local correction to the solution provided by the pre-computed reduced model. We will see in the examples that, as the conjugate gradient does not need to search for the part of the solution that is directly captured in the pre-computed reduced space, the convergence rate is very high.

The diagonal preconditioner $\precond^{-1}$ classically accounts for local heterogeneities in the structure. A much more sophisticated preconditioner can of course be used (e.g.:incomplete Cholesky factorization). Yet, if the reduced basis is correctly pre-computed and/or updated, the computational effort required to obtain a good solution of the projected problem should be relatively low. Hence, the potentially expensive construction of an advanced preconditioner is not \textit{a priori} justified. The example section will show that the number of conjugate gradient iterations required to compute a correction to the reduced model is indeed kept very low, even while only using this crude preconditioner.


\section{Application to the adaptive reduction of damage problems}
\label{sec:application_damage}


\subsection{Damageable lattice model}
\label{sec:lattice}

\begin{figure}[p]
       \centering
       \includegraphics[width=0.5 \linewidth]{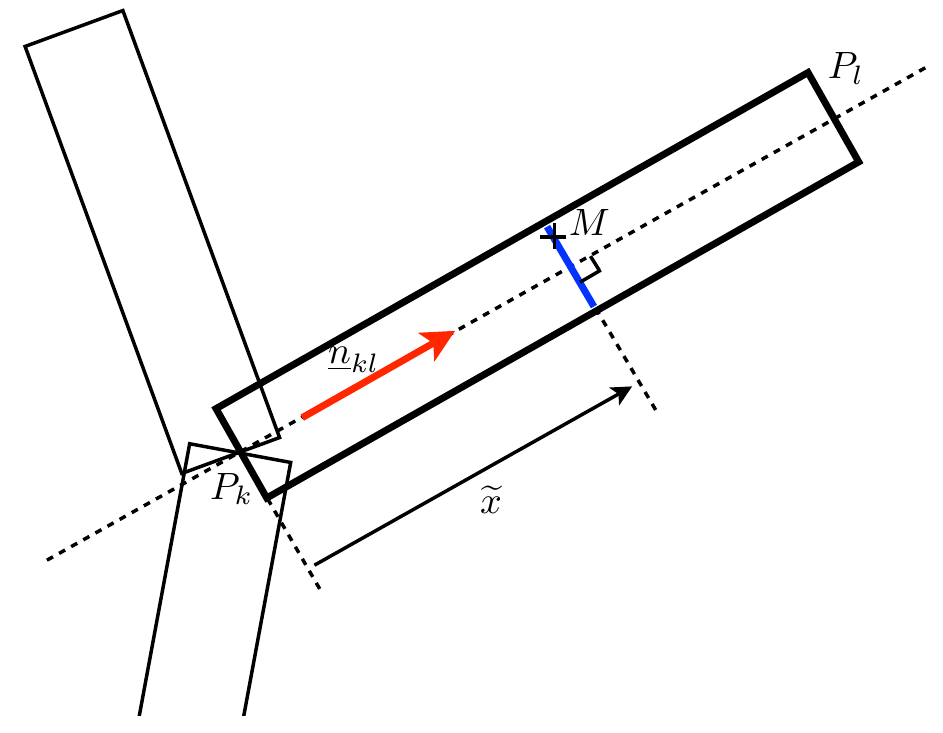}
       \caption{Definition of the lattice problem.}
       \label{fig:lattice}
\end{figure}

\begin{figure}[p]
       \centering
       \includegraphics[width=0.9 \linewidth]{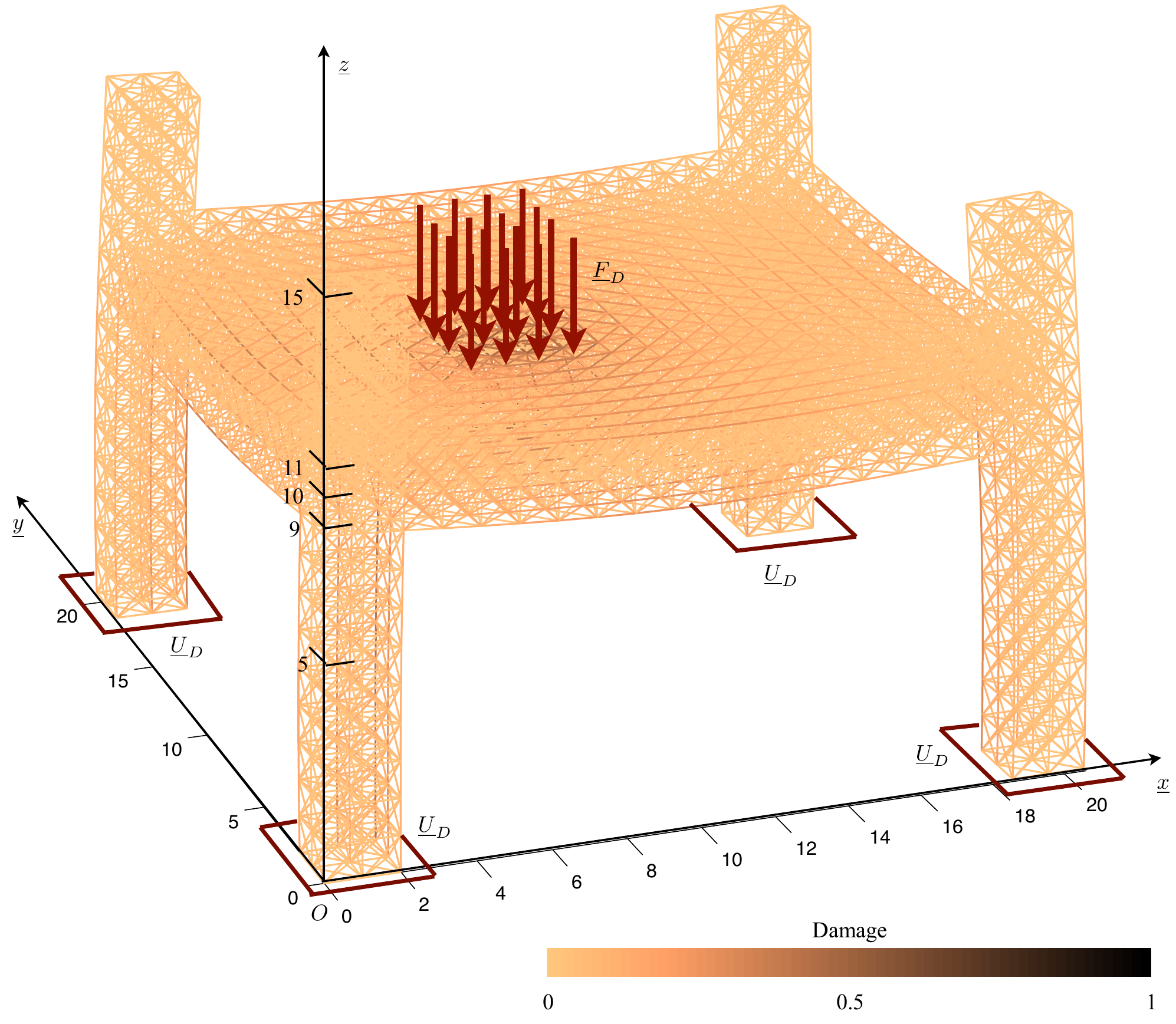}
       \caption{Reference solution obtained at the eighth load increment of the full-scale simulation (maximum value of prescribed forces before unstable damage propagation). Darker bars have a higher damage state.}
       \label{fig:Reference}
\end{figure}

We focus on a very simple damageable lattice structure, made of bars under traction or compression (figure (\ref{fig:Reference})). Elementary bar indiced $b$ occupies domain $\Omega_{b}$ such that
$\Omega = \bigcup_{b \in \llbracket 1 , n_b\rrbracket } \Omega_b$, where $n_b$ is the number of bars. The direction of bar number $b$ is denoted by $\V{n}_{kl}=\V{P_k P_l}/\|\V{P_k P_l}\|$ where $P_k$ and $P_l$ ($k<l$ and $(k,l)\in \llbracket 1 , n_P\rrbracket^2$) are the two extremities of the bar as defined in figure (\ref{fig:lattice}). The displacement field in is supposed constant in a cross-section. Therefore, the displacement of any point $M \in \Omega_b$ can be expressed in the following form: 
\begin{equation}
 \disp_{|M\in \Omega_b} = \left( \disp(\tilde{x}) . \V{n}_{kl} \right) \V{n}_{kl}  = u(\tilde{x}) \, \V{n}_{kl} 
\end{equation}
where each cross section is parametrised by a scalar $\tilde{x} \in [0 , \lVert \V{P_k P_l} \rVert]$, defined by $\tilde{x} =\V{ P_kM } . \V{n}_{kl}$. 

We use a constitutive law based on classical damage mechanics to describe the behaviour of the lattice structure. The lineic strain energy of bar $b$ reads:
\begin{equation}
\displaystyle {e_d}_{|\tilde{x}} = \frac{ 1 }{ 2 } E(1-d)  \, S_b \, \epsilon_{|\tilde{x}}^2
\end{equation}
where $\epsilon_{|\tilde{x}} = u_{,\tilde{x}}(\tilde{x})$ is the strain measured in the direction of the bar, $(S_b)_{b \in \llbracket 1 , n_b\rrbracket }$, is the section of the bar, and $d$ is a damage variable which ranges from 0 (undamaged material) to 1 
(completely damaged material). The local state equations, derived from the strain energy read:
\begin{equation}
\left\{ \begin{array}{l}
\displaystyle N_{|\tilde{x}} = \frac{ \partial e_d }{ \partial \epsilon } = E(1-d)  \, S_b  \, \epsilon_{|\tilde{x}} 
\\
\displaystyle  Y_{|\tilde{x}}  = -\frac{ \partial e_d }{ \partial d } =  \frac{ 1 }{ 2 } E  \, S_b  \, \epsilon_{|\tilde{x}}^2
\end{array} \right.
\end{equation}
where we have introduced the axial force $N = S_b \left( \V{n}_{kl}^{T} . \stress .\V{n}_{kl} \right)$, and the thermodynamic force $Y$.
A non-local evolution law is defined to link the damage variables to the thermodynamic forces:
\begin{equation}
d_{|\tilde{x}} = \max_{\tilde{x} \in [0 , \lVert \V{P_k P_l} \rVert]} \left( \min \left( 1,\max_{\tau \leq t} \left( \alpha  ( Y_{|\tilde{x},\tau})^\beta \right) \right) \right)
\end{equation}

We suppose that the lattice is not subjected to any volume force, that the material of each elementary bar is isotropic and homogeneous, and that the section of each bar is constant. As a consequence, the resulting displacement field evolves linearly along the direction of each bar.

In terms of finite element discretization, each bar will be considered as a linear element
\begin{equation}
\disp(\tilde{x}) = \left( 1 - \frac{\tilde{x}}{\|\V{P_k P_l}\|} \right) \disp_{|P_l} + \frac{\tilde{x}}{\|\V{P_k P_l}\|} \, \disp_{|P_k} 
\end{equation}
Because of the assumptions made previously on the loading, geometry and material properties, the finite element solution yields the ``exact'' solution to the damageable lattice problem. \\

\subsection{Lattice structure and loading}

We solve the problem defined in figures (\ref{fig:Reference}). The sections and Young's modulus are unitary. 
The material parameters are $\beta=0.5$ and $\alpha=\sqrt{2}$. Any bar parallel to axis $(O,\V{x})$, $(O,\V{y})$ or $(O,\V{z})$ has a unitary length. 

The Dirichlet boundary conditions applied to the structure read, for any node $(P_i)$ (with $i \in \llbracket 1, n_p\rrbracket$) of the lattice such that $(P_i)\in (O,\V{x},\V{y})$:
\begin{equation}
\disp_{|P_i } . \V{z} = 0
\end{equation}
As represented in figure (\ref{fig:Reference}), homogenous Neuman boundary conditions, in the $-\V{z}$ direction, are applied to any point $P_i$ that belongs to plane $\mathcal{P}$ defined by:
\begin{equation}
\mathcal{P} = \{ M = (x,y,z) \, | \, x \in [7 , 9] \textrm{ and } y \in [8 , 10] \textrm{ and } z=11 \}
\end{equation}

\begin{figure}[htb]
       \centering
       \includegraphics[width=0.7 \linewidth]{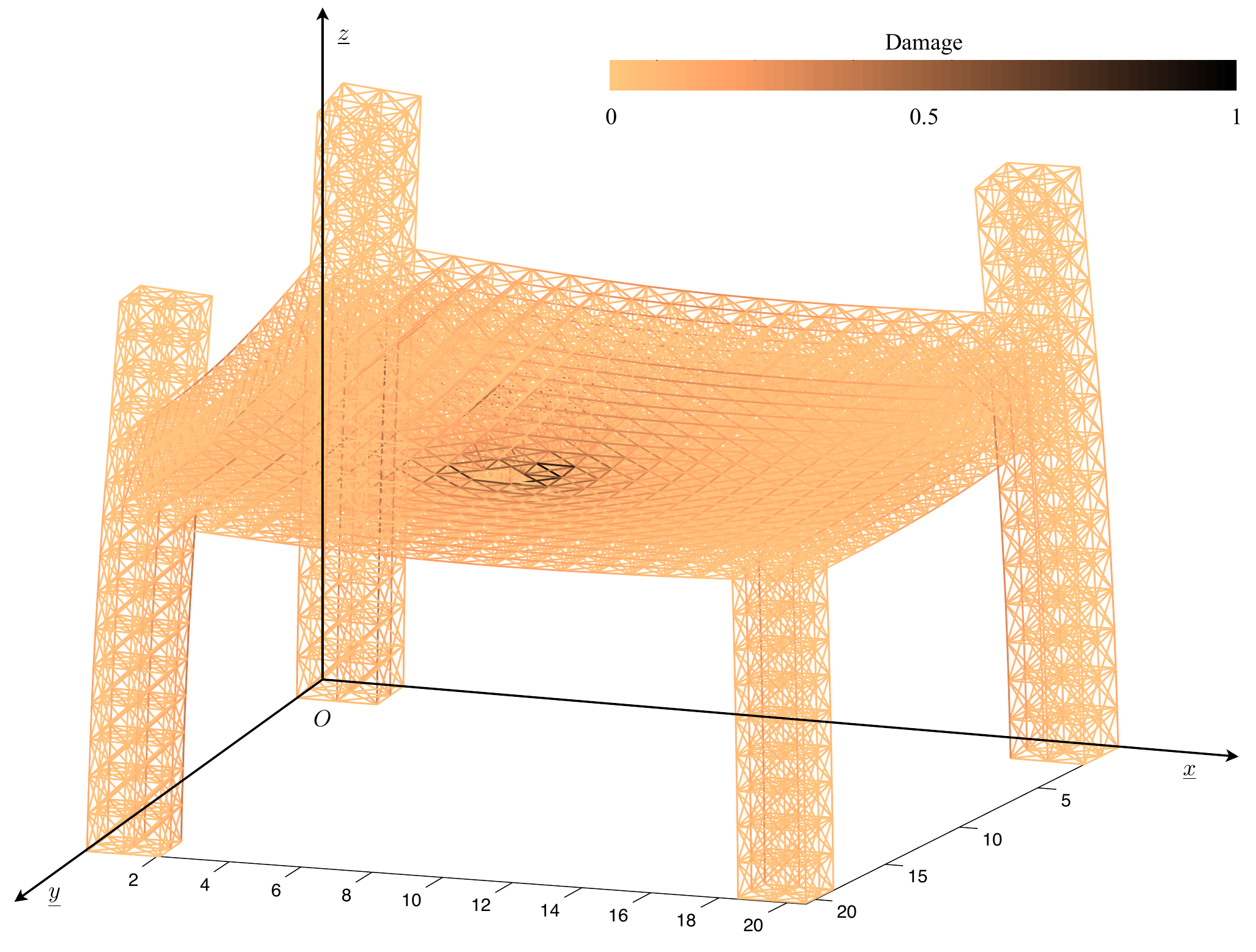}
       \caption{Final damage obained at the end of the full-scale simulation, performed in 30 load increments. The bars that are completely damaged are not represented.}
       \label{fig:Reference2}
\end{figure}

Increasing the value of the applied forces leads to the initiation and propagation of a damaged zone as shown in figure (\ref{fig:Reference2}). This phenomenon is unstable from a structural point of view, and the damage state represented in this figure is past the limit point at which the force applied to the structure reaches its maximum value. An arc-length algorithm with local control is used to solve the damage propagation beyond this limit point. The basic idea of this algorithm is to look for the amplitude of the force under the constraint that the maximum increment in the damage variables over a time step is fixed (see \cite{riks1972,allixcorigliano1996} or \cite{kerfridengosselet2010} for additional details). The reference solution is obtained by solving the full problem \eqref{eq:discretized_equilibrium} in 30 of these damage increments of the Newton/arc-length algorithm. The final state is represented in figure (\ref{fig:Reference2}). 

Our purpose is to check the adequacy of the proposed local/global reduction technique when trying to obtain at cheap costs the load/deflection curve (pre and post failure phases) of the structure under the loading conditions described previously. The results provided in the following sections give the value of the norm of the external forces vector (also called loading factor) obtained by the arc-length method. The displacement used to plot the load/deflection curves is the average vertical displacement of the points at which forces are applied.

\subsection{Particularisation of the local/global approach}

\subsubsection{Definition of the Local/Global splitting}
\label{sec:splitting}

\begin{figure}[htb]
       \centering
       \includegraphics[width=1.0 \linewidth]{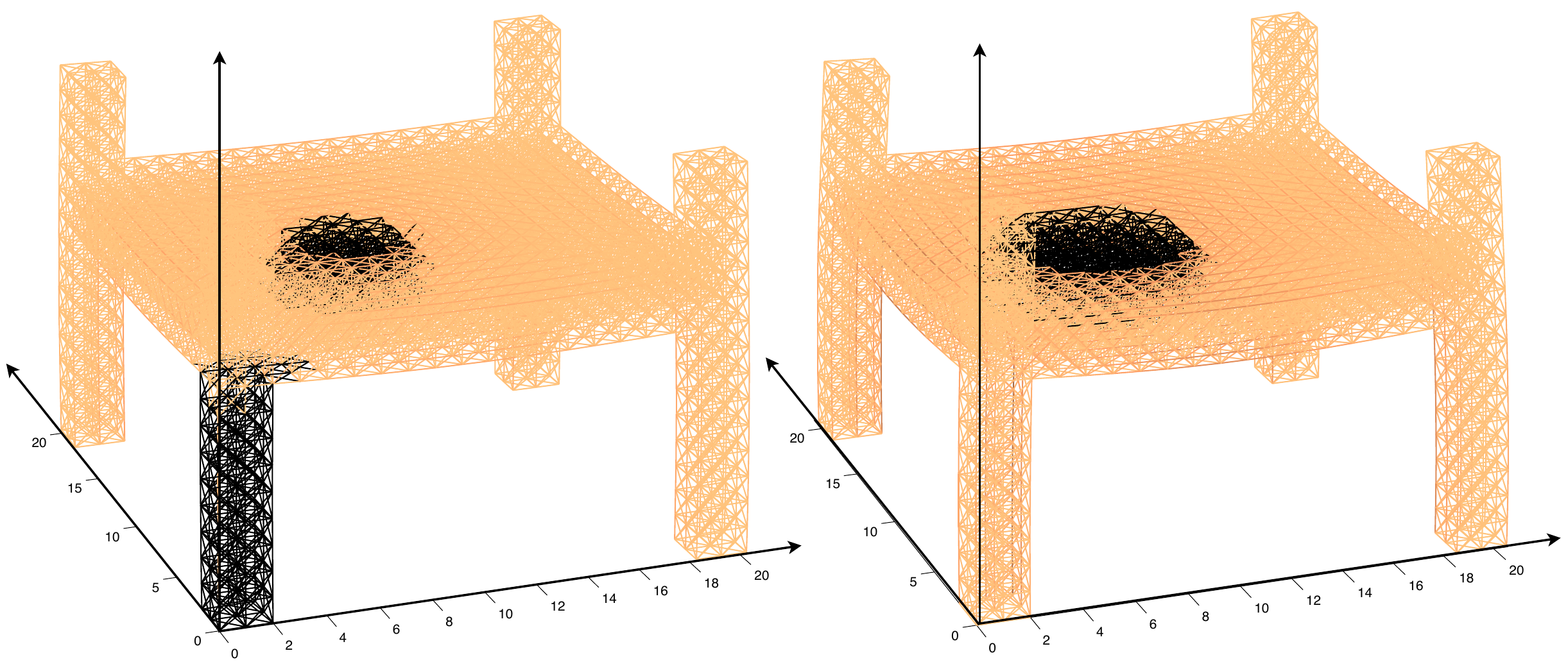}
       \caption{Splitting obtained at two different stages of the simulation performed by the local/global reduction method (left: load increment 4, right: load increment 30). Dark bars are connected to at least one fully resolved node. As time increases, the degrees of freedom for which no reduction is performed localises in the region of failure.}
       \label{fig:Patches}
\end{figure}

We need at that point to particularize the general method developed in section \ref{sec:LocGlob}. The splitting between reduced degrees of freedom $\widetilde{\Deltadispdiscr}^{(r)}$ and the complementary ones $\widetilde{\Deltadispdiscr}^{(f)}$ is made by considering that the part of the structure undergoing strong damage variations will not be correctly solved for when projected on a pre-computed global basis (see \cite{gendreallix2009,kerfridenallix2009} for similar assumptions in a domain decomposition-based local/global approaches). This statement is purely based on mechanical intuition and will be validated in the following examples.

The following procedure is adopted. The splitting will be fixed over a time step, to avoid the stagnation or divergence of the Newton process. At the end of a time step, the element undergoing the maximum damage increment is spotted, and a sphere of radius $\rho_s$, centred at the isobarycenter of the element is created. Every degree of freedom belonging to a node located inside this sphere is set as a fully resolved degree of freedom.

This procedure is repeated on the remaining elements (those which do not have all their degrees of freedom in the previous list).

The algorithm is stopped if one of the two following statements is satisfied.
\begin{itemize}
\item The maximum damage increment in the remaining elements is lower than $k_\textrm{Dam}$ times the global maximum. In our examples, $k_\textrm{Dam}$ is set to $0.5$. This ensures that the fully resolved degrees of freedom are actually connected to elements undergoing a significant increase in their damage state compared to the remaining of the structure.
\item The number of fully resolved degrees of freedom exceeds $k_\textrm{LocGlo}$ times the total number of unknowns. $k_\textrm{LocGlo}$ is set to 0.1 in our examples. This statement will be satisfied if the former one is not, which means that the structure undergoes a virtually homogeneous increase in its damage state. In that case, the local/global procedure might be inefficient.
\end{itemize}

To illustrate this procedure, two different reduced domains are represented in figure (\ref{fig:Patches}). The dark elements are connected to at least one fully resolved node. The first one is obtained at the end of the third time step and is used to solve the fourth time step of the analysis, while the second one is used to solve the thirtieth and last time step. The first figure is obtained during the initiation phase. Notice that the maximum damage increments are found in the zone where the the load is applied, and in the ``pillar'' that is closer to this zone, which undergoes a significant bending loading. In the later case, corresponding to unstable propagation phase, the damage increment localises in the zone where the load is applied.

This empirical choice for the local/global splitting is only justified by its efficiency. However we only considered one type of model for brittle fracture. We expect the optimal splitting to depend on the chosen model and on the stage of the structural failure. A more suitable criterion should be investigated in future work, based on algebra to circumvent this dependency. For instance, one could consider the decrease rate of the singular values of the correlation matrix restricted to an assumed smoothed domain and loosely maximise it, with a greedy-like algorithm, by successive small modifications of this domain.

\subsubsection{Modification of the reduced basis}
\label{sec:global_update}

At the end of a time step, the reduced basis is enriched in order to take into account the information provided by the local iterative solver. The procedure used here is simply to add the successive solutions to the initial snapshot, and sort them by a singular value computation. This step can be optimized, following our previous work \cite{kerfridengosselet2010}. The issue here is that the reduced basis is only locally enhanced, which might results in ill-posed reduced problems if not correctly addressed. This question will not be further discussed in this paper, which focuses on the benefits provided by the local/global approach.

\subsection{\textit{A posteriori} reduced basis}

\subsubsection{Reference solution to construct the \textit{a posteriori} reduced basis}

\begin{figure}[htb]
       \centering
       \includegraphics[width=0.6 \linewidth]{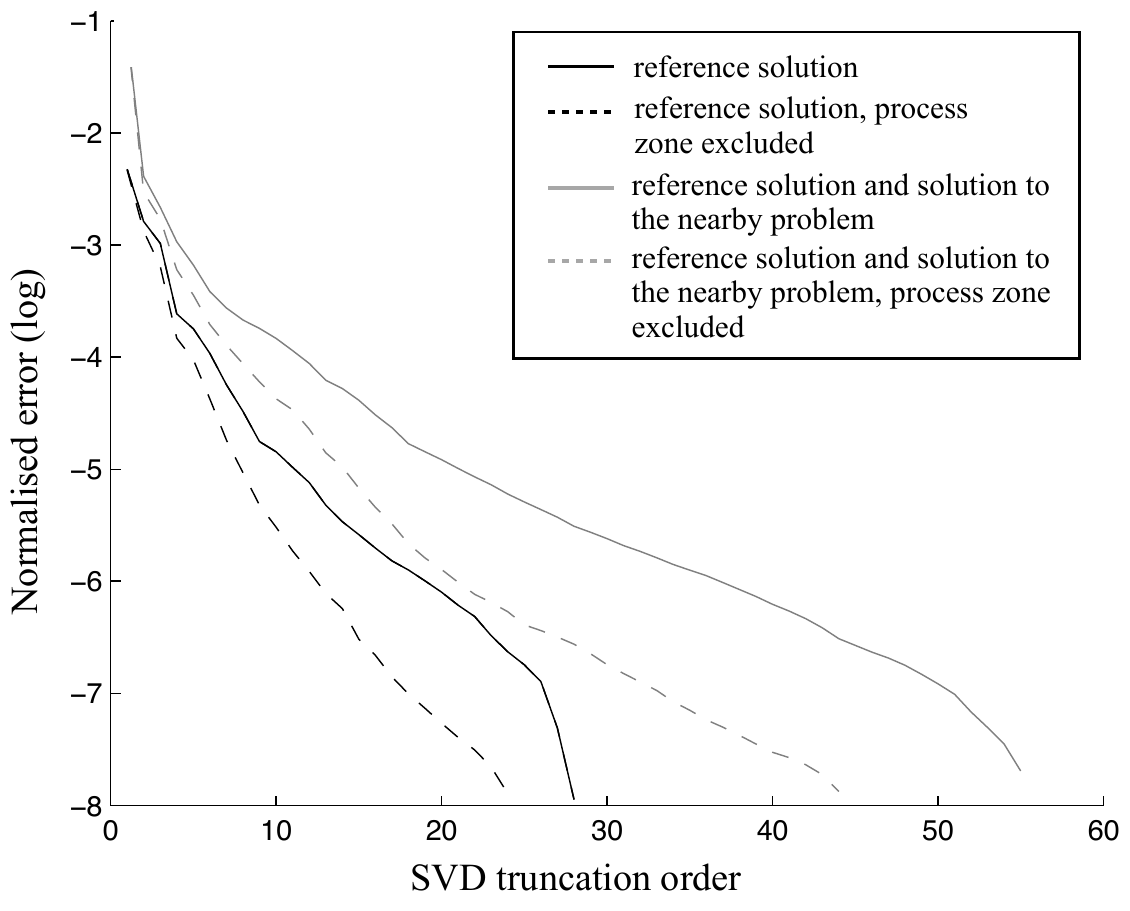}
       \caption{
Measure of the error in the snapshot representation as a function of the error of truncation of the singular value decomposition. In solid lines dark line, the set of 30 reference solutions are defined as the snapshot. In grey lines, a second snapshot operator is defined by concatenating the $30$ reference solutions and the $30$ solutions of the nearby problem. Both analyses are repeated when excluding the process zone (defined \textit{a priori} by the set of points $\{ M = (x,y,z) \, | \, x \in [5 , 13] \textrm{ and } y \in [6 , 13] \})$) from the snapshot definition (i.e.: the snapshot now is the restriction of the solutions to the part of the domain which does not undergo significant damage). The resulting error curves are plotted in dashed lines.
}
       \label{fig:SVD}
\end{figure}

\begin{figure}[htb]
       \centering
       \includegraphics[width=1.0 \linewidth]{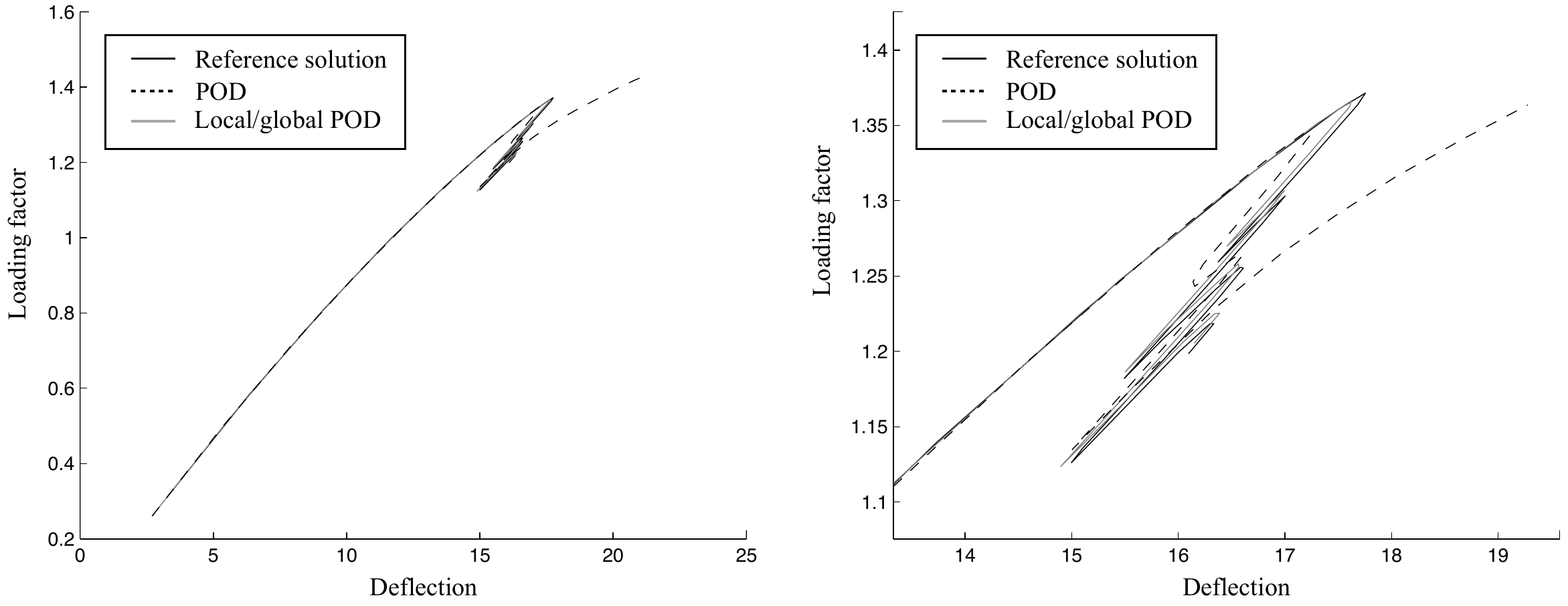}
       \caption{Load/deflection curves obtained one using successively direct numerical simulation, simple proper orthogonal decomposition and local/global reduction technique. The snapshot used here is the solution of the direct numerical simulation.}
       \label{fig:Curve_POD_LocGlob1}
\end{figure}

The first test consists in using as a snapshot the consecutive reference solutions obtained by a Newton/arc-length solution algorithm performed on the initial system of equations, without reduction. In other terms, we try to reproduce the snapshot by reduced order modelling. It terms of realistic applications, this example is not of particular interest as one would need to know the solution in advance in order to obtain the followinf results. However, studying this academic case is interesting to understand the behaviour of the proposed strategy. 

The snapshot is composed of $n_s = 30$ vectors. This information is compressed by using a singular value decomposition as detailed in section \ref{sec:POD}. We plot, in figure (\ref{fig:SVD}), the normalised error in the snapshot representation
\begin{equation}
\nu_{\textrm{SVD}} =  \frac{ \displaystyle \left(  \sum_{j=1}^{n_s} \left\Vert  \V{\mathbf{S}}^{j} - \sum_{i=1}^{n_c} \left( {\V{\mathbf{C}}^{i}}^T \V{\mathbf{S}}^{j} \right) \V{\mathbf{C}}^i \right\Vert^2 \right)^{\frac{1}{2}} }{ \displaystyle \left( \sum_{j=1}^{n_s} \left\Vert  \V{\mathbf{S}}^{j}  \right\Vert^2 \right)^{\frac{1}{2}} }  = \frac{ \displaystyle \left( \sum_{i=n_c+1}^{n_s} \lag_i  \right)^{\frac{1}{2}} }{ \displaystyle \left(  \sum_{i=1}^{n_s} \lag_i  \right)^{\frac{1}{2}} } 
\end{equation}
as a function of the order of truncation of the SVD, where the eigenvalues $(\lambda_I)_{i \in \llbracket 1, n_s\rrbracket}$ of the correlation matrix are now sorted in decreasing order. The decrease rate of this criterion is relatively low, which suggests that the information from the snapshot is poorly suited to compression. For reference, in \cite{ryckelynck2005}, the author estimates that the snapshot space is correctly approximated by the space spanned by the first $n_c$ singular vectors if $\nu_{\textrm{SVD}} \leq 10^{-8}$. In our case, this would lead to the selection of the $30$ singular vectors to define a relevant reduced model. The singular vector basis is here truncated at order $n_c = 3$ in order to yield a significant (thus observable) error in the solution when using the basic POD method. This error can then be compared to the one obtained when using the proposed local/global scheme. \\

In figure (\ref{fig:Curve_POD_LocGlob1}), the reference load/deflection curved is compared to the one obtained when using the basic POD on the first hand, and when using the local/global reduction technique on the other hand. The peak load obtained with the basic POD is slightly underestimated, by less than 2\%, but the overall behaviour of the structure is correctly captured. The local/global algorithm permits to further increase the accuracy of the approximated solution. The error in the maximum load factor is now less than 0.5 \%. \\

\begin{figure}[htb]
       \centering
       \includegraphics[width=1 \linewidth]{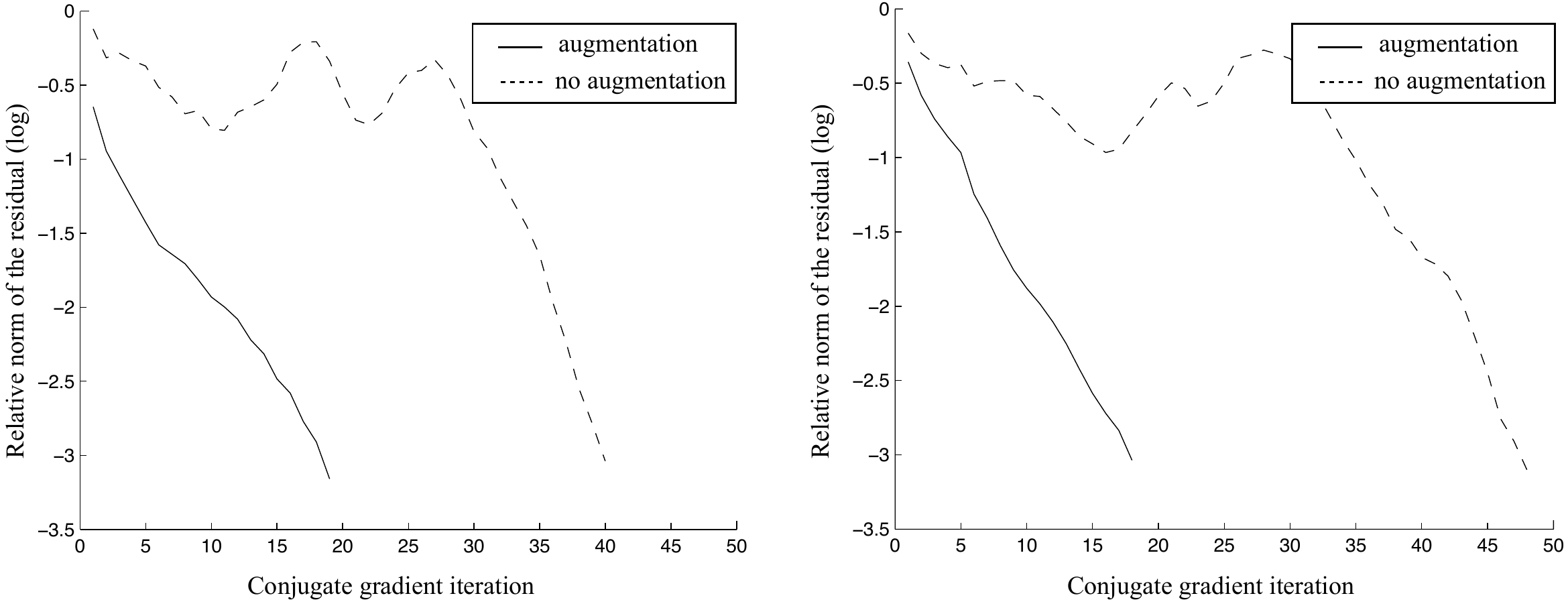}
       \caption{Convergence curves of the augmented conjugate gradient used to solve the local problem, at the third Newton iteration of load increment 15. On the left-hand side, the reference solution is used as a snapshot. On the right-hand side, the snapshot used to define the reduced model is the solution of the nearby problem.}
       \label{fig:CG}
\end{figure}

We now analyse the effect of preconditioning the local iterations of the conjugate gradient by the projection technique detailed in section \ref{sec:iterativeCG}. The convergence curves that are plotted in figure (\ref{fig:CG}, left) are obtained at the third Newton iteration of the fiftieth time step of the analysis (i.e.: during the damage propagation phase). The norm of the condensed system \eqref{eq:condensed} of equation, normalized by the norm of the condensed right-hand term is observed as a function of the number of iterations of the conjugate gradient solver. One can see that the number of iteration to convergence is significantly decreased when using the augmented Krylov solver.

An other interesting observation is that the convergence curve decreases in a monotonic manner, which means that the error in the global residual keeps decreasing as a function of the cost of the local conjugate gradient. This is not, of course, a proof, and this conclusion is probably problem dependent, but the same behaviour as been observed in all our numerical experiments. A probable reason for this is that the augmentation deflates the spectrum of the condensed operator, removing its highest eigenvalues and forcing the Krylov algorithm to its super-convergence phase (one can refer to \cite{gosselet2003} for a proper analysis of the effect of eigenvectors-based augmentation on the convergence of Krylov solvers).

\subsubsection{Nearby solution to construct the \textit{a posteriori} reduced basis}
\label{sec:snap2}

\begin{figure}[p]
       \centering
       \includegraphics[width=1.0 \linewidth]{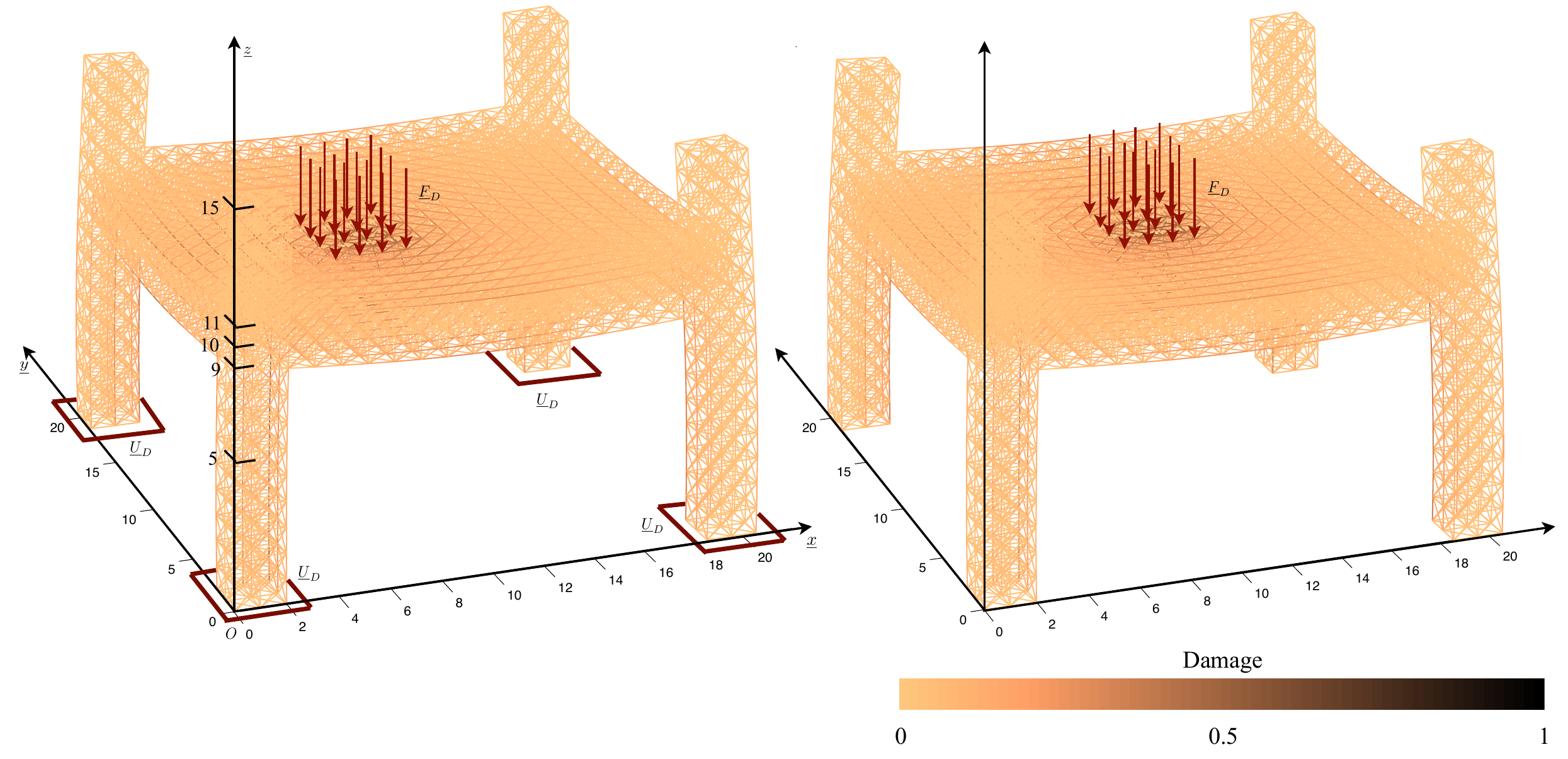}
       \caption{Two set of reference solutions obtained by direct numerical simulations and used as snapshot vectors in our examples. Only the solutions to the eighth load increment are represented in both cases. On the left-hand side, the snapshot is the reference solution which is looked for. On the right-hand side, the snapshot is the solution to a slightly different problem (``nearby problem''), where the position of the applied load has been modified.}
       \label{fig:Snap}
\end{figure}

\begin{figure}[p]
       \centering
       \includegraphics[width=0.6 \linewidth]{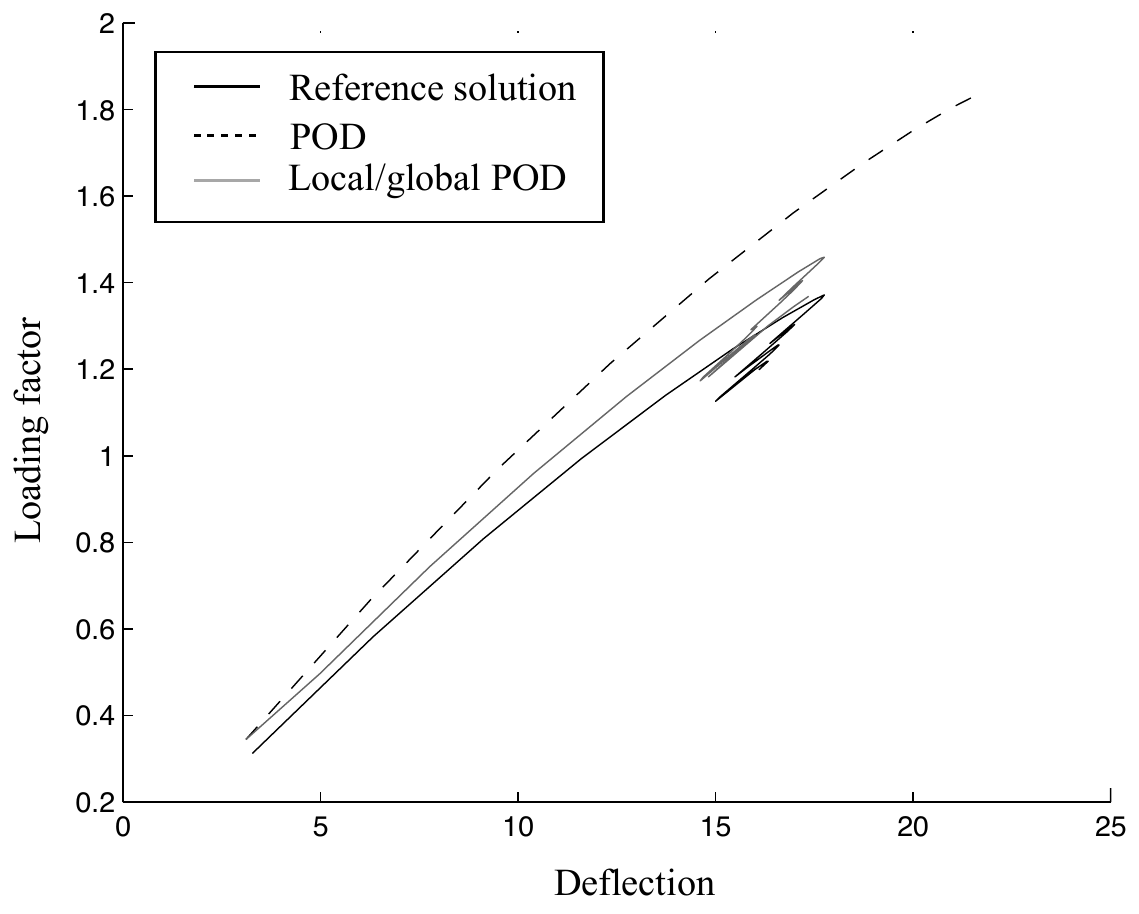}
       \caption{Load/deflection curves obtained one using successively direct numerical simulation, simple proper orthogonal decomposition and local/global reduction technique. The snapshot used here is the solution of the direct numerical simulation to the nearby problem.}
       \label{fig:Curve_POD_LocGlob2}
\end{figure}

In the second test case, we use as a snapshot the consecutive solutions obtained by solving a nearby problem, which differs from the reference one only by the position of the applied forces: homogenous Neumann boundary conditions are applied to any point $P_i$ (with $i \in \llbracket 1, n_P\rrbracket$) such that $P_i \in\mathcal{P}'$ defined by:
\begin{equation}
\mathcal{P}' = \{ M = (x,y,z) \, | \, x \in [9 , 11] \textrm{ and } y \in [9 , 11] \textrm{ and } z=11 \}
\end{equation}

This snapshot is for instance a particular output of a parametric study meant to obtain the response of the system subjected to various load cases. One would expect to be able to reuse the information generated during this simulation to fasten the solution process of the next one. The difference between the snapshot solution and the reference solution that is looked for is illustrated in figure (\ref{fig:Snap}). It should be noted that far away from the part of the structure where Neumann boundary conditions are applied, the solutions are qualitatively similar. Yet, they are very different within this process zone, which qualitatively justifies the local/global splitting proposed in section \eqref{sec:splitting}. The same effect appears in figure (\ref{fig:SVD}, which shows that the combination of the information from the two load cases can be significantly compressed when removing the process zone from the SVD analysis.

Figure (\ref{fig:Curve_POD_LocGlob2}) compares the load/deflection curves obtained in the reference case, and when using successively the basic POD strategy and the local/global reduction technique. Obviously, the POD strategy looses its relevance very quickly, as the damage localises in the part of the structure where the damage would localise during the snapshot simulation. The resulting model is too stiff, and the maximum strength of the structure is overestimated. When using the local/global approach, with the parameters described previsouly for the fine analysis of the damaged areas, a much more accurate load/deflection curve is obtained. Yet the peak load is still overestimated by 5\%. This remaining inaccuracy is due to the fact that the infomation from the snapshot is not correct in the reduced zone. Global corrections of the reduced basis are still necessary. This issue is discussed in the next section. \\

In figure (\ref{fig:CG}, right) we analyse the effect of preconditioning the local iterations of the conjugate gradient by the projection technique. As in the previous test case, the convergence curves are observed at the third Newton Newton iteration of the fiftieth time step of the analysis. Again, the number of iteration to convergence is significantly decreased when using the augmentation. We would expect this precontionning technique to be less efficient in this particular case. Yet it is not the case, which is probably due to the fact that the global update described in section \ref{sec:global_update} is mechanically relevant.

\section{Adaptive reduced basis}

We now try to reduce the inaccuracy observed in figure (\ref{fig:Curve_POD_LocGlob2}) by using global corrections of the reduced model. The adaptation algorithm used here was introduced in \cite{kerfridengosselet2010}. We just recall the basics of this particular strategy, and explain how it can be used to enhance the proposed local/global scheme.

\subsection{"On-the-fly" global corrections of the reduced basis}
\label{sec:global_correc}

The principle of these global corrections, is to update the reduced basis used to define reduced problem \eqref{eq:discr_problem_reduc}. It is, for now, not linked to the global/local technique described previously in this paper. The updates are done ``on-the-fly", at any iteration $i+1$ of the Newton algorithm used to solve \eqref{eq:discr_problem_reduc} at a given time step $t_{n+1}$ of the analysis. If the reduced problem is sufficiently converged, and if the residual of the norm of the initial problem evaluated at iteration $i$ is estimated too high, then the following linearized system is considered:
\begin{equation}
\label{eq:linearized_global_correc}
\optangbar^i \bar{\deltadispdiscr} = - \resdiscr^{i}
\end{equation}
$\optangbar^{i}$ is an approximation of the tangent $\optang^{i}$, while $\resdiscr^{i}$ is the residual of the initial system of equations \eqref{eq:discretized_equilibrium} evaluated in $\dispdiscr_{|t_n} + \basereduc^i \, \coeffreduc^i$. Operator $\basereduc^i$ designates the reduced basis considered at iteration $i$ of the Newton process.

The solution to this problem is searched for in two supplementary spaces by a projected Krylov algorithm:
\begin{equation}
\begin{array}{l}
\displaystyle \bar{\deltadispdiscr}= \bar{\deltadispdiscr}_\textrm{C} + \bar{\deltadispdiscr}_\textrm{K} \\
\displaystyle \textrm{where} \quad
\left\{ \begin{array}{l}
\displaystyle \bar{\deltadispdiscr}_\textrm{C} \in \textrm{Im}(\basereduc^i) \\
\displaystyle \bar{\deltadispdiscr}_\textrm{K} \in \textrm{Im} (\basereduc^i)^{\bot_{\optangbar^i}} = \textrm{Ker}({\basereduc^i}^T \optangbar^i)
\end{array} \right.
\end{array}
\end{equation}
This problem is solved coarsely by a projected Krylov solver ($10^{-1}$ is the typical value of the stopping criterion), othogonally to the current reduced basis. The resulting solution, which is $\optangbar^i-$orthogonal to the $\basereduc^i$, is added to the reduced basis:
\begin{equation}
\basereduc^{i+1} = \left( \basereduc^{i} \quad \frac{\bar{\deltadispdiscr}_\textrm{K}}{ \lVert \bar{\deltadispdiscr}_\textrm{K} \rVert} \right)
\end{equation}
Newton prediction $i+1$ can now be solved after updating residual $\resdiscrred^i$ and reduced stiffness $\optangred^i$ with respect to this new basis vector.

Various features are added to this algorithm to efficiently sort the successively added basis vectors, but the main idea described above is sufficient to understand the following results.

\subsection{Coupling with the local/global approach}

Once the global norm of the residual of locally reduced problem \eqref{eq:discr_problem_reduc_loc_glo} is sufficiently converged, and if the global residual is juged too high, a correction vector to the global reduced basis is calculated, following the algorithm given in the previous subsection. The residual and tangent operator of problem \eqref{eq:discr_problem_reduc_loc_glo} are updated and condensed tangent problem \eqref{eq:condensed} is solved with the projected conjugate gradient algorithm. This algorithm converges in a very few iterations as the new basis vector, added to the augmentation space, is very close to the solution which is looked for. The slight mismatch is only due to:
\begin{itemize}
\item the approximation of the tangent used to provide the global correction vector.
\item the coarse solution of \eqref{eq:linearized_global_correc}.
\item the approximation of the displacement in the reduced zone of the local/global approach.
\end{itemize}

\subsection{Results}

\begin{figure}[p]
       \centering
       \includegraphics[width=1.0 \linewidth]{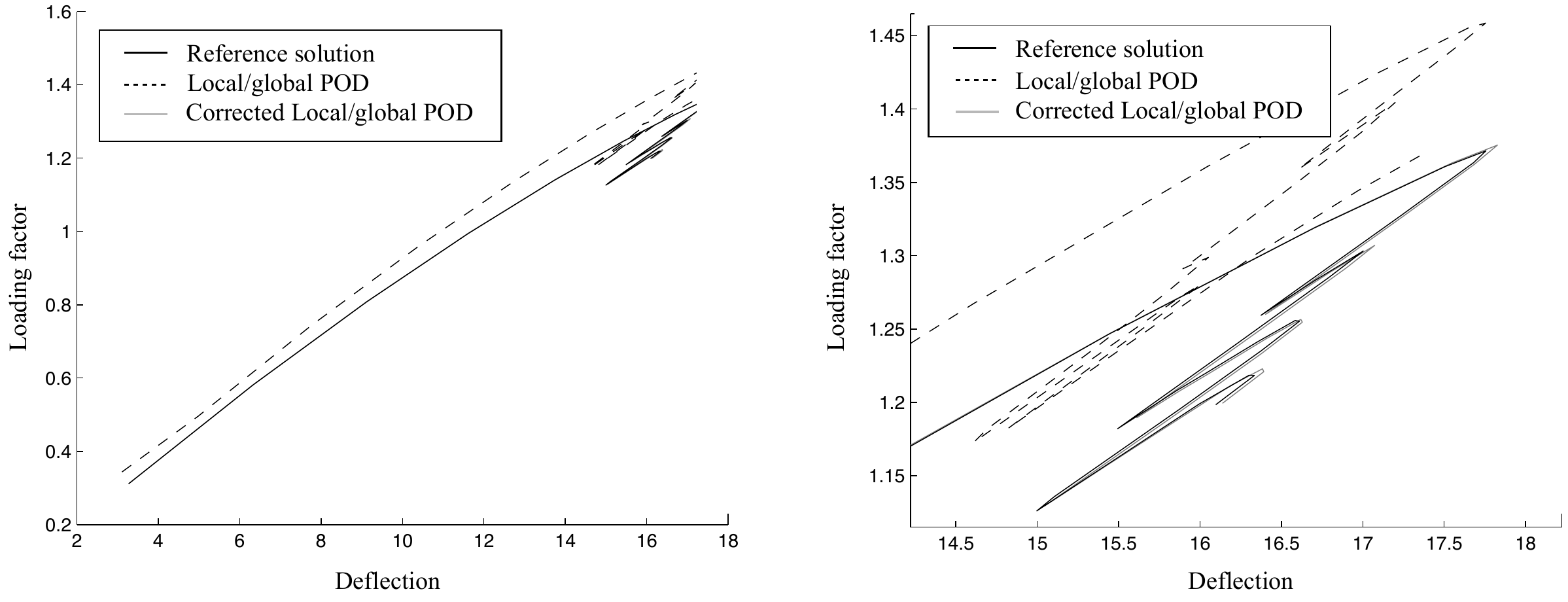}
       \caption{Load/deflection curves obtained one using successively direct numerical simulation, local/global reduction technique and local/global reduction technique with ''on-the-fly" global corrections.}
       \label{fig:Curve_POD_LocGlob_CG}
\end{figure}

\begin{figure}[p]
       \centering
       \includegraphics[width=0.6 \linewidth]{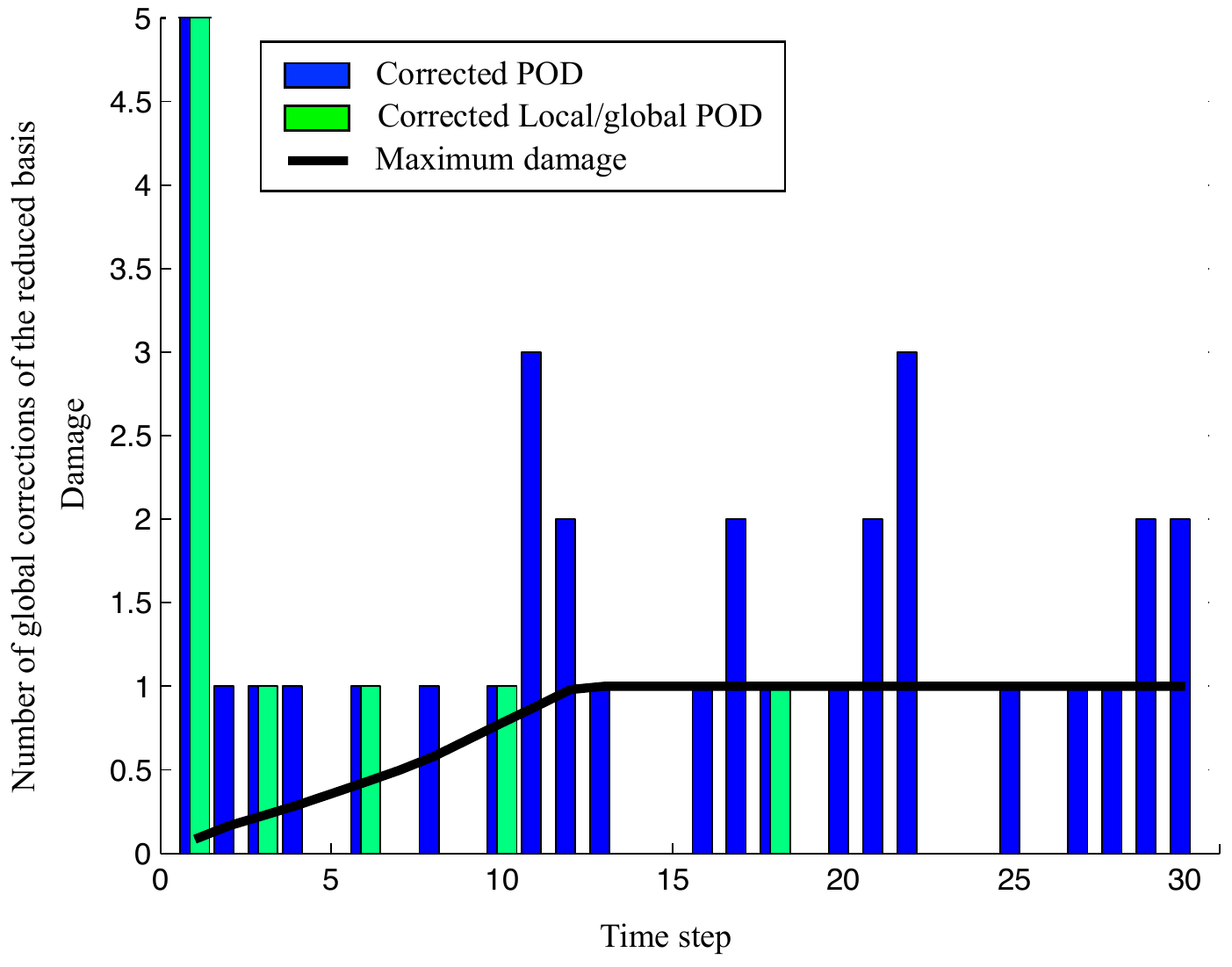}
       \caption{Number of global corrections performed when using the adaptive model order reduction on the one hand, and when using adapative local/global model order reduction on the other hand.}
       \label{fig:IterGlo}
\end{figure}

\begin{figure}[p]
       \centering
       \includegraphics[width=1. \linewidth]{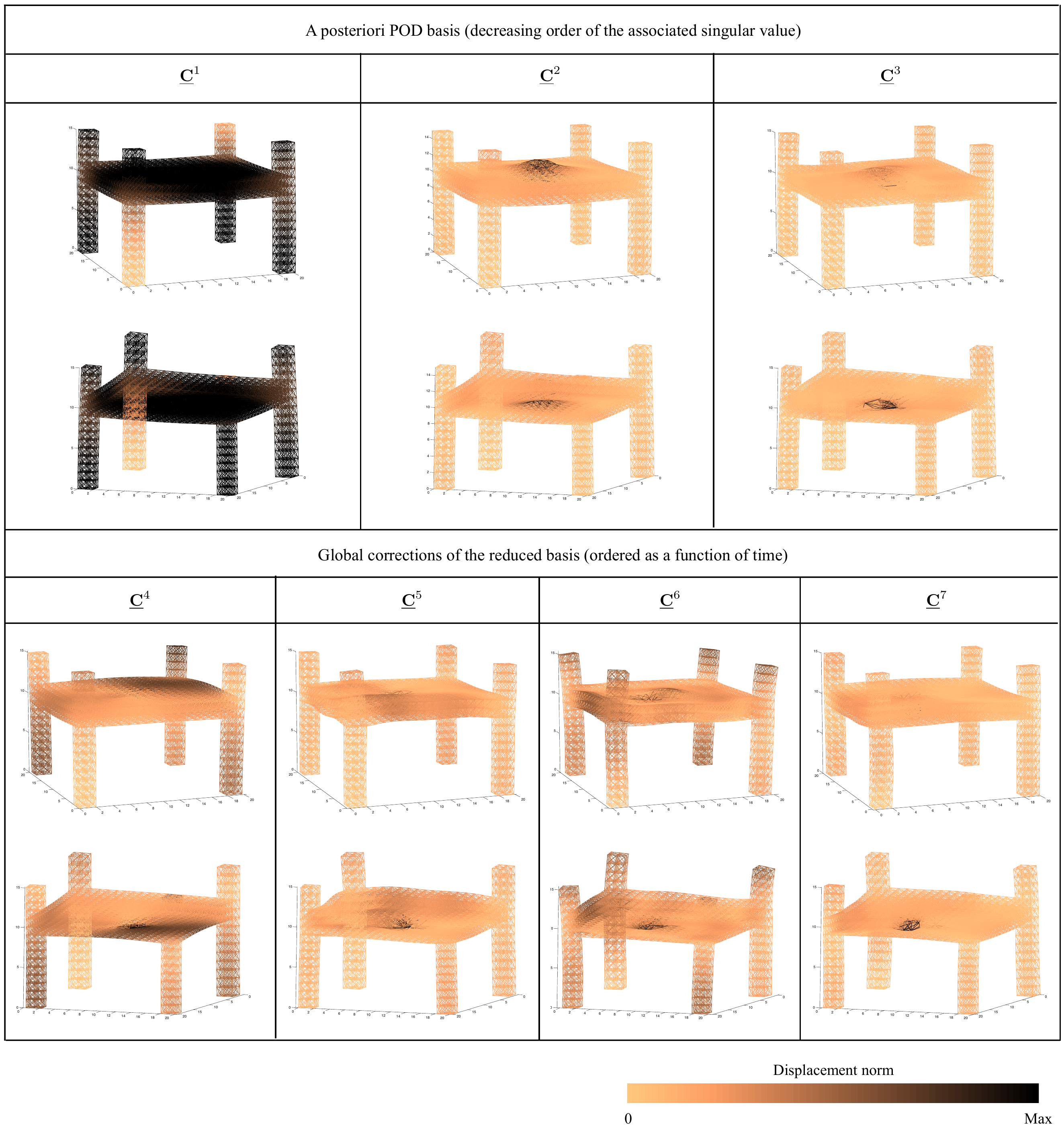}
       \caption{A posteriori POD basis (3 vectors) and successive global corrections (4 vectors computed respectively at time steps 1, 3, 10 and 18). Two different views of each basis vector are represented for better qualitative understanding.}
       \label{fig:CorrectedBasis}
\end{figure}

The  test case described in section \ref{sec:snap2} is solved using successively the local/global technique on its own, the global correction algorithm on its own, and the local/global reduction technique coupled with global corrections. The target norm of the relative global residual norm is $10^{-1}$. Corrections are allowed only if the relative norm of the residual of locally reduced problem \eqref{eq:discr_problem_reduc_loc_glo} is below $10^{-3}$. It is clear from the load-displacement curve in Figure 11 that the global corrections significantly improve the results of the local/global algorithm.  (\ref{fig:Curve_POD_LocGlob_CG}). An illustration of the POD basis and its successive corrections is given in figure (\ref{fig:CorrectedBasis}).

More importantly, the number of relatively expensive global corrections required to achieve the target in terms of global residual norm is very low when using the local/global approach, as opposed to the one observed when using only global corrections. This effect is shown in figure (\ref{fig:IterGlo}). On can also notice that using only the global correction methods lead to an increase in the computational costs in the damage propagation phase of the analysis, which has been reported in \cite{kerfridengosselet2010}. With the local/global reduction approach, this effect disappears. The same trends are observed in \cite{kerfridenallix2009}, where a local/global algorithm was developed in a domain decomposition framework to avoid unnecessary computations far away from the process zones. This suggests that the result obtained here may be extended, and that the number of global corrections required to achieve a given level of accuracy does not depend on the local nonlinearity (i.e: independent on the damage state).

\section{Conclusions and discussion}

The work described in this paper focuses on decreasing the computational effort required in solving large scale damage and fracture problems, where small scale phenomena must be taken into account to accurately represent the global structural behaviour. Although most investigations on this issue rely on the extension of homogenization frameworks, we followed an alternative route, choosing projection-based model order reduction as a starting point.
The proposed local/global reduction technique was demonstrated on a damaging 3D frame structure, but we expect that the algorithm will carry forward to the more general case of the initiation and propagation of cracks in structures.


The main conclusions of the paper are:
\begin{itemize}
\item The proposed local/global model order reduction approach can significantly improve the relevancy of using a global reduced basis to approximate the displacement field in the case of localised failure;
\item The number of global corrections of the reduced model required to obtain a given level of accuracy is drastically reduced when the balance equations which exhibit the highest level of nonlinearity are excluded from the reduction. 
\end{itemize}

In most of the applications for which this methodology is intended, the additional computational costs are kept minimal for two reasons:
\begin{itemize}
\item The damaged zone where no reduction is performed is small compared to the size of the structure in realistic engineering problems;
\item The local fine scale solution is efficiently obtained by Krylov corrections of a good initialisation provided by the restriction of the reduced model to the damaged zone. 
\end{itemize}

Yet, some issues need to be addressed before global/local reduction algorithms can be efficiently and safely used to tackle engineering problem:
\begin{itemize}
\item Cost of the construction of the reduced systems. A system reduction is needed to ensure the computational efficiency of the method. It is now well known that projection-based model order reduction only reduces the size of the algebraic systems to solve, and that the construction of these systems is still expensive. Indeed, it involves the evaluation of the internal forces to construct the tangent stiffness and the residual of the reduced problem, which requires global updates of the constitutive behaviour, and integrations over the whole domain. System approximations permit to solve this issue. They basically consist in obtaining a cheap, but reliable approximation of the internal forces. In \cite{kerfridengosselet2010}, we used the so-called hyperreduction \cite{ryckelynck2005} (also introduced as the missing-point approximation by  \cite{astridweiland2008} in a different context). We found that adaptive reduced basis methods need to be carefully tailored when used within this framework. Indeed, one cannot rely on the fact that the residual of the full-scale problem is known at any iteration of the nonlinear solver anymore: the number of global evaluation of the internal forces must be kept  minimal. In addition, we observed that the stability of the particular system approximation that we used is highly dependent on the choice of its parameters (reduced integration domain), especially when the damage mechanism localises. We foresee that excluding the most damaged zones from the reduced domain can significantly improve the stability of this scheme. Some numerical investigations are required to confirm this intuitive speculation.
\item Quality control. Each step in the algorithm is a source of error: (i) choice of the set of degrees of freedom for which the solution is searched in a reduced space (ii) choice of the reduced basis (iii) choice of the accuracy for the local solution (iv) choice of the system reduction (v) choice of the time stepping parameters.
These choices should be made within an encompassing framework to ensure that the error associated with each of them can be controlled in the same manner. Therefore, one needs to find a relevant error criterion and build a control process, capable of ensuring that the above choices introduce the same order of error in the solution so that the global error be not dominated by any of them in particular.
\end{itemize} 

In terms of further applications, this method should be validated using more advanced models for brittle fracture, in particular nonlocal continuum damage or softening plasticity models, for which the distinction between local and global effects is not trivial. The extension to ductile failure and/or geometric nonlinearities should also be possible, the global corrections being expected to efficiently handle the ``smooth'' part of the nonlinearities, while the local/global scheme should treat separately the localisation of defects.

We believe that this local/global reduction framework is an important milestone to obtain an optimal reduction method applicable to the parallel simulation of failure in solids. Indeed, as alluded to in the introduction of this paper, several authors have proposed to perform a systematic model order reduction of the expensive local problems (so-called ``problems by substructure'') classically involved in domain decomposition methods \cite{rixen2004,markovicibrahimbegovic2009,ladevezepassieux2009}. The application of such ideas to failure assessment has never been tackled so far, as it first requires an efficient reduction scheme for the local problems, such as the one proposed in this paper, in order to guarantee reasonable scalabilities and load balancing.

\section{Acknowledgements}

Pierre Kerfriden would like to thank Pierre Gosselet (LMT/ENS Cachan, France) for his precious advices. St\'ephane Bordas and Pierre Kerfriden would like to acknowledge the financial support of the Royal Academy of Engineering and of the Leverhulme Trust for Bordas' Senior Research Fellowship ``Towards the next generation surgical simulators'' (2009-2010) as well as the support of EPSRC under grant EP/G042705/1 Increased Reliability for Industrially Relevant Automatic Crack Growth Simulation with the eXtended Finite Element Method.
The support of the British Council Exchange Programme is also gratefully acknowledged by the authors.

\bibliographystyle{unsrt}
\bibliography{bibliography}

\end{document}